# Studies in Tours of Knight on Rectangular Boards

Awani Kumar (Lucknow, India)

E-mail: awanieva@gmail.com

**Abstract**

*The author has constructed and enumerated tours of knight having various magic properties on 4 x n and 6 x n boards. 16 magic tours of knight have been discovered on 4 x 18 board, 88 on 4 x 20 board, 464 on 4 x 22 board, 2076 on 4 x 24 board, 9904 on 4 x 26 board and 47456 on 4 x 28 board. Magic tours exist on all boards of size 4 x 2k for k > 8. Quasi-magic tour exists on 6 x 11 board. 8 magic tours of knight have been discovered on 6 x 12 board and magic tours exist on all boards of size 6 x 4k for k > 2.*

**Introduction:** British puzzlist H. E. Dudeney [1] said, "The knight is the irresponsible low comedian of the chessboard." Knight is a unique piece and its crooked move has enchanted and fascinated not only chess players but also attracted attention of mathematics enthusiasts. It has been the basis for so many unusual and amusing combinatorial problems. "The oldest of knight puzzles", asserts mathematics popularizer M. Gardner [2] "is the knight's tour". Books, namely, by Rouse Ball [3], Kraitchik [4], Pickover [5], Petkovic [6], Wells [7] and by many others cover tour of knight problem. Some great mathematical minds, such as, Raimond de Montmort, Abraham de Moivre, and Leonhard Euler have delved into the knight's tour problem, but mostly on 8 x 8 board. According to Jelliss [8], "oldest knight tour on 8x8 board is by Ali C. Mani, an otherwise unknown chess player, and al-Adli ar-Rumi who flourished around 840 AD and is known to have written a book on Shatranj (the form of chess then popular)". Warnsdorf [9] proposed an interesting method – "Play the knight to a square where it commands the fewest cells not yet used" and Stonebridge [10] improved it further. Jelliss [11] has also described methods given by Roget and Collini but all these methods can't be used to enumerate or estimate the number of tours. Recently Parberry [12] has described an algorithm for constructing closed knight's tour on square boards. Kyek et al [13] have described an algorithm to estimate the bounds of the number of closed tours on a square board. Ismail et al [14] have given a new algorithm for constructing knight's tour on 8 x 8 board. Although Cull [15] and McGown [16] have given proofs and examples of knight's tours on rectangular boards and Baskoro [17] has also given a method for constructing closed knight's tour on rectangular boards but most of the available literature on knight's tour is related to square boards. Rectangular boards have got scanty attention. Construction of simple knight's tours on 4 by n (n > 4) boards is relatively easy and the number of knight's tours on smaller boards (up to 4 x 10) has been counted. But, what about knight tours on larger rectangular boards? Can we have magic knight's tours on 4 by n boards? If yes, what is the smallest board size and how many magic tours are there? It is the tours having magic properties and the enumeration of tours, which have attracted, enchanted and fascinated innumerable people and the author plans to look into them.

## A. Tours of knight on 4 by n boards

**Description:** Like its conventional meaning, knight tour on a 4 x n board is consecutive knight moves where knight visits all the 4n cells without visiting any cell twice. If the initial cell of the tour and its last cell are connected by knight's move, it is called re-entrant or closed knight tour else it is open knight tour. Generally, knight tours on larger board have fewer restrictions but, surprisingly, knight tours on 4 x n board have more restrictions than the 3 x n board. Sainte-Maire [18] proved that re-entrant knight tours are impossible on 4 x n board. However, closed knight tours are possible on 3 x n board if n is even and n >= 10. Similarly, there are 4 distinct knight tours (3 in geometrical form) on 3 x 4 board but no knight tour is possible on 4 x 4 board. Before proceeding further, it is better to explain few terms associated with knight tours. A '**magic knight tour**' on a rectangular board 4 x n (n even) has all the 4 long rows adding to magic constant n(4n+1)/2 and the n short rows adding to magic constant 2(4n+1). There can't be magic knight tour on a rectangular board 4 x n if n is odd because the adjacent long rows will contain an odd and an even number of odd-numbered cells. So the sum of the rows will be alternately odd and even. A '**semi-magic knight tour**' has either all the long rows adding up to a constant sum, called magic constant (MC) or all the short rows adding up to a different magic constant. Jelliss [19] has defined two special types of semi-magic tours, namely, '**quasi-magic**' tour in which the non-magic lines add to only two different values and '**near-magic**' tour in which the non-magic lines add to the magic constant and two other values. Such types of semi-magic tours are of interest here because magic knight tours are not there on smaller boards. Sum of short and long rows have been shown in the figures and magic constants have been highlighted in bold font.

**4 x 9 board**: Boards up to 4 x 8 size have already been studied in sufficient details. Jelliss [20] mentions, "There are 112 symmetric open tours". Figure 1 shows a semi-magic tour with short rows adding to magic constant 74.

| 1 | 18 | 35 | 20 | **74** | | 1 | 4 | 19 | 36 | 60 | | 1 | 18 | 35 | 20 | **74** |
|---|---|---|---|---|---|---|---|---|---|---|---|---|---|---|---|---|
| 34 | 21 | 16 | 3 | **74** | | 32 | 35 | 2 | 5 | **74** | | 34 | 21 | 2 | 17 | **74** |
| 17 | 2 | 19 | 36 | **74** | | 3 | 18 | 33 | 20 | **74** | | 9 | 16 | 19 | 36 | 80 |
| 22 | 33 | 4 | 15 | **74** | | 34 | 31 | 6 | 17 | 88 | | 22 | 33 | 10 | 3 | 68 |
| 9 | 14 | 23 | 28 | **74** | | 7 | 16 | 21 | 30 | **74** | | 15 | 8 | 23 | 28 | **74** |
| 32 | 27 | 10 | 5 | **74** | | 22 | 29 | 8 | 15 | **74** | | 32 | 27 | 4 | 11 | **74** |
| 13 | 8 | 29 | 24 | **74** | | 11 | 14 | 23 | 26 | **74** | | 7 | 14 | 29 | 24 | **74** |
| 26 | 31 | 6 | 11 | **74** | | 28 | 25 | 12 | 9 | **74** | | 26 | 31 | 12 | 5 | **74** |
| 7 | 12 | 25 | 30 | **74** | | 13 | 10 | 27 | 24 | **74** | | 13 | 6 | 25 | 30 | **74** |
| 161 | 166 | 167 | 172 | | | 151 | 182 | 151 | 182 | | | 159 | 174 | 159 | 174 | |

Fig.1. Semi-magic tours     Fig.2. Tours with long rows adding to alternate values.

Discerning readers must have seen that long rows in Figure 1 are adding to four different values. This is common in tours on 4 x n rectangular boards with odd side. Figure 2 shows two tours with long (9-cell) rows adding up to only two different values. There are 1682 such tours.

**4 x 10 board**: Regarding semi-magic tours, Jelliss [19] asserts that "This case has not been looked at so far." The author has enumerated 3102 semi-magic tours, 30 quasi-magic tours and 22 near-magic tours with long (10-cell) rows magic. They have been shown in Figure 3.

| 7 | 10 | 23 | 30 | 70 | 5 | 8 | 25 | 28 | 66 | 3 | 18 | 37 | 24 | **82** |
|---|----|----|----|-----|----|----|----|----|-----|----|----|----|----|----|
| 22 | 29 | 8 | 11 | 70 | 26 | 29 | 4 | 7 | 66 | 38 | 23 | 2 | 19 | **82** |
| 9 | 6 | 31 | 24 | 70 | 9 | 6 | 27 | 24 | 66 | 17 | 4 | 25 | 36 | **82** |
| 28 | 21 | 12 | 19 | 80 | 30 | 23 | 10 | 3 | 66 | 22 | 39 | 20 | 1 | **82** |
| 5 | 18 | 25 | 32 | 80 | 11 | 2 | 21 | 32 | 66 | 5 | 16 | 35 | 26 | **82** |
| 40 | 27 | 20 | 13 | 100 | 22 | 31 | 12 | 1 | 66 | 40 | 21 | 6 | 15 | **82** |
| 17 | 4 | 33 | 26 | 80 | 13 | 20 | 33 | 40 | 106 | 7 | 14 | 27 | 34 | **82** |
| 36 | 39 | 14 | 1 | 90 | 34 | 37 | 16 | 19 | 106 | 28 | 31 | 8 | 11 | 78 |
| 3 | 16 | 37 | 34 | 90 | 17 | 14 | 39 | 36 | 106 | 13 | 10 | 33 | 30 | 86 |
| 38 | 35 | 2 | 15 | 90 | 38 | 35 | 18 | 15 | 106 | 32 | 29 | 12 | 9 | **82** |
| **205** | **205** | **205** | **205** | | **205** | **205** | **205** | **205** | | **205** | **205** | **205** | **205** | |

Fig.3. Semi-magic, quasi-magic and near-magic tours with long rows magic.

There are 104 semi-magic tours, 12 quasi-magic tours and 16 near-magic tours with short (4-cell) rows magic. They have been shown in Figure 4.

| 7 | 10 | 31 | 34 | **82** | 3 | 18 | 39 | 22 | **82** | 1 | 20 | 39 | 22 | **82** |
|---|----|----|----|----|----|----|----|----|----|----|----|----|----|----|
| 30 | 33 | 8 | 11 | **82** | 40 | 21 | 2 | 19 | **82** | 38 | 23 | 18 | 3 | **82** |
| 9 | 6 | 35 | 32 | **82** | 17 | 4 | 23 | 38 | **82** | 19 | 2 | 21 | 40 | **82** |
| 26 | 29 | 12 | 15 | **82** | 24 | 37 | 20 | 1 | **82** | 24 | 37 | 4 | 17 | **82** |
| 5 | 14 | 27 | 36 | **82** | 5 | 16 | 35 | 26 | **82** | 5 | 16 | 35 | 26 | **82** |
| 28 | 25 | 16 | 13 | **82** | 36 | 25 | 6 | 15 | **82** | 36 | 25 | 6 | 15 | **82** |
| 17 | 4 | 37 | 24 | **82** | 7 | 14 | 27 | 34 | **82** | 7 | 14 | 27 | 34 | **82** |
| 22 | 39 | 20 | 1 | **82** | 28 | 31 | 10 | 13 | **82** | 28 | 31 | 10 | 13 | **82** |
| 3 | 18 | 23 | 38 | **82** | 11 | 8 | 33 | 30 | **82** | 11 | 8 | 33 | 30 | **82** |
| 40 | 21 | 2 | 19 | **82** | 32 | 29 | 12 | 9 | **82** | 32 | 29 | 12 | 9 | **82** |
| 187 | 199 | 211 | 223 | | 203 | 203 | 207 | 207 | | 201 | **205** | **205** | 209 | |

Fig.4. Semi-magic, quasi-magic and near-magic tours with short rows magic.

**4 x 11 board**: Jelliss [20] mentions that the "number of tours (on 4 x 11 board is) not determined." The author has computed 3341926 arithmetically distinct tours. Out of these tours, there are 267 semi-magic tours with short rows magic as shown in Figure 5. Figure 6 shows two interesting tours with the long rows adding up to only two values (247 and 248 here). There are 57586 such tours.

| 1 | 22 | 43 | 24 | **90** | | 1 | 20 | 43 | 24 | 88 | | 3 | 20 | 43 | 24 | **90** |
|---|---|---|---|---|---|---|---|---|---|---|---|---|---|---|---|---|
| 42 | 25 | 20 | 3 | **90** | | 44 | 23 | 2 | 21 | **90** | | 42 | 23 | 2 | 21 | 88 |
| 21 | 2 | 23 | 44 | **90** | | 19 | 4 | 25 | 42 | **90** | | 19 | 4 | 25 | 44 | 92 |
| 26 | 41 | 4 | 19 | **90** | | 26 | 41 | 22 | 3 | 92 | | 26 | 41 | 22 | 1 | **90** |
| 5 | 18 | 27 | 40 | **90** | | 5 | 18 | 27 | 40 | **90** | | 5 | 18 | 27 | 40 | **90** |
| 36 | 39 | 6 | 9 | **90** | | 36 | 39 | 6 | 9 | **90** | | 36 | 39 | 6 | 9 | **90** |
| 17 | 8 | 37 | 28 | **90** | | 17 | 8 | 37 | 28 | **90** | | 17 | 8 | 37 | 28 | **90** |
| 38 | 35 | 10 | 7 | **90** | | 38 | 35 | 10 | 7 | **90** | | 38 | 35 | 10 | 7 | **90** |
| 13 | 16 | 29 | 32 | **90** | | 13 | 16 | 29 | 32 | **90** | | 13 | 16 | 29 | 32 | **90** |
| 34 | 31 | 14 | 11 | **90** | | 34 | 31 | 14 | 11 | **90** | | 34 | 31 | 14 | 11 | **90** |
| 15 | 12 | 33 | 30 | **90** | | 15 | 12 | 33 | 30 | **90** | | 15 | 12 | 33 | 30 | **90** |
| 248 | 249 | 246 | 247 | | | 248 | 247 | 248 | 247 | | | 248 | 247 | 248 | 247 | |

Fig.5. Semi-magic tour    Fig.6. Tours with long rows adding up to alternate values.

**4 x 12 board**: In search for magic-tour, Jelliss [19] has examined only tours of squares and diamonds type and found 16 semi-magic tours (4 of them being quasi-magic) with short (4-cell) rows magic. Jean-Charles Meyrignac computed all the 48 quasi-magic tours with long (12-cell) rows magic in 2003. The author has computed 18243164 arithmetically distinct tours. There are 608 semi-magic, 28 quasi-magic and 72 near-magic tours with short rows magic. There are 58356 semi-magic, 48 quasi-magic and 112 near-magic tours with long rows magic. Near-magic tours haven't been investigated earlier. Figure 7 shows three such tours.

| 1 | 24 | 47 | 26 | **98** | | 1 | 24 | 47 | 26 | **98** | | 11 | 8 | 41 | 38 | **98** |
|---|---|---|---|---|---|---|---|---|---|---|---|---|---|---|---|---|
| 46 | 27 | 2 | 23 | **98** | | 46 | 27 | 2 | 23 | **98** | | 42 | 39 | 10 | 7 | **98** |
| 21 | 4 | 25 | 48 | **98** | | 21 | 4 | 25 | 48 | **98** | | 9 | 12 | 37 | 40 | **98** |
| 28 | 45 | 22 | 3 | **98** | | 28 | 45 | 22 | 3 | **98** | | 34 | 43 | 6 | 15 | **98** |
| 5 | 20 | 43 | 30 | **98** | | 5 | 20 | 43 | 30 | **98** | | 13 | 16 | 33 | 36 | **98** |
| 44 | 29 | 18 | 7 | **98** | | 44 | 29 | 18 | 7 | **98** | | 44 | 35 | 14 | 5 | **98** |
| 19 | 6 | 31 | 42 | **98** | | 19 | 6 | 31 | 42 | **98** | | 17 | 2 | 47 | 32 | **98** |
| 32 | 41 | 8 | 17 | **98** | | 32 | 41 | 8 | 17 | **98** | | 48 | 45 | 4 | 1 | **98** |
| 9 | 16 | 33 | 40 | **98** | | 9 | 16 | 33 | 40 | **98** | | 3 | 18 | 31 | 46 | **98** |
| 34 | 37 | 12 | 15 | **98** | | 36 | 39 | 10 | 13 | **98** | | 26 | 29 | 24 | 21 | 100 |
| 13 | 10 | 39 | 36 | **98** | | 15 | 12 | 37 | 34 | **98** | | 19 | 22 | 27 | 30 | **98** |
| 38 | 35 | 14 | 11 | **98** | | 38 | 35 | 14 | 11 | **98** | | 28 | 25 | 20 | 23 | 96 |
| 290 | **294** | **294** | 298 | | | **294** | 298 | 290 | **294** | | | **294** | **294** | **294** | **294** | |

Fig.7. Near-magic tours on 4 x 12 board.

**4 x 13 board**: Tours on 4x13 board haven't been investigated earlier. The author has computed 100641235 arithmetically distinct tours. There are 1444 semi-magic tours with short rows magic as shown in Figure 8. Here, the long rows, summing up to four consecutive numbers (namely, 343 to 346), have an aesthetic appeal. There are no quasi-magic tours with short rows magic.

| 1 | 26 | 51 | 28 | **106** |
|---|---|---|---|---|
| 50 | 29 | 24 | 3 | **106** |
| 25 | 2 | 27 | 52 | **106** |
| 30 | 49 | 4 | 23 | **106** |
| 5 | 22 | 47 | 32 | **106** |
| 48 | 31 | 6 | 21 | **106** |
| 7 | 20 | 33 | 46 | **106** |
| 42 | 45 | 8 | 11 | **106** |
| 19 | 10 | 43 | 34 | **106** |
| 44 | 41 | 12 | 9 | **106** |
| 15 | 18 | 35 | 38 | **106** |
| 40 | 37 | 16 | 13 | **106** |
| 17 | 14 | 39 | 36 | **106** |

343 344 345 346

| 1 | 26 | 51 | 28 | **106** |
|---|---|---|---|---|
| 50 | 29 | 24 | 3 | **106** |
| 25 | 2 | 27 | 52 | **106** |
| 30 | 49 | 4 | 23 | **106** |
| 21 | 6 | 31 | 48 | **106** |
| 32 | 47 | 22 | 5 | **106** |
| 7 | 20 | 33 | 46 | **106** |
| 42 | 45 | 8 | 11 | **106** |
| 19 | 10 | 43 | 34 | **106** |
| 44 | 41 | 12 | 9 | **106** |
| 15 | 18 | 35 | 38 | **106** |
| 40 | 37 | 16 | 13 | **106** |
| 17 | 14 | 39 | 36 | **106** |

343 344 345 346

| 1 | 26 | 51 | 28 | 106 |
|---|---|---|---|---|
| 52 | 29 | 24 | 3 | 108 |
| 25 | 2 | 27 | 50 | 104 |
| 30 | 49 | 4 | 23 | **106** |
| 5 | 22 | 47 | 32 | **106** |
| 48 | 31 | 6 | 21 | **106** |
| 7 | 20 | 33 | 46 | **106** |
| 42 | 45 | 8 | 11 | **106** |
| 19 | 10 | 43 | 34 | **106** |
| 44 | 41 | 12 | 9 | **106** |
| 15 | 18 | 35 | 38 | **106** |
| 40 | 37 | 16 | 13 | **106** |
| 17 | 14 | 39 | 36 | **106** |

345 344 345 344

Fig.8. Semi-magic tours with short rows magic.   Fig.9

There are 1563998 tours with long rows adding up to only two different values. Figure 9 is an example of such tours (with long rows adding to 344 and 345). Discerning readers must have noted that the long rows are also adding up to consecutive numbers. There are 122848 such tours.

**4 x 14 board**: Tours on 4 x 14 board haven't been investigated earlier. The author has computed 526152992 arithmetically distinct tours. There are 3480 semi-magic tours, 136 quasi-magic tours and 244 near-magic tours with short rows magic. They are shown in Figure 10 to Figure 12.

| 1 | 54 | 3 | 32 | 5 | 50 | 9 | 34 | 21 | 38 | 11 | 40 | 15 | 44 | 357 |
|---|---|---|---|---|---|---|---|---|---|---|---|---|---|---|
| 28 | 31 | 26 | 53 | 8 | 33 | 6 | 37 | 10 | 35 | 18 | 43 | 12 | 41 | 381 |
| 55 | 2 | 29 | 4 | 49 | 24 | 51 | 20 | 47 | 22 | 39 | 14 | 45 | 16 | 417 |
| 30 | 27 | 56 | 25 | 52 | 7 | 48 | 23 | 36 | 19 | 46 | 17 | 42 | 13 | 441 |
| **114** | **114** | **114** | **114** | **114** | **114** | **114** | **114** | **114** | **114** | **114** | **114** | **114** | **114** | |

Fig.10. Semi-magic tour on 4 x 14 board with short rows adding up to magic constant 114.

| 1 | 54 | 25 | 32 | 5 | 52 | 23 | 48 | 9 | 38 | 11 | 40 | 15 | 44 | 397 |
|---|---|---|---|---|---|---|---|---|---|---|---|---|---|---|
| 28 | 31 | 4 | 53 | 24 | 33 | 6 | 37 | 20 | 47 | 18 | 43 | 12 | 41 | 397 |
| 55 | 2 | 29 | 26 | 51 | 22 | 35 | 8 | 49 | 10 | 39 | 14 | 45 | 16 | 401 |
| 30 | 27 | 56 | 3 | 34 | 7 | 50 | 21 | 36 | 19 | 46 | 17 | 42 | 13 | 401 |
| **114** | **114** | **114** | **114** | **114** | **114** | **114** | **114** | **114** | **114** | **114** | **114** | **114** | **114** | |

Fig.11. Quasi-magic tour on 4 x 14 board with short rows adding to magic constant 114.

| 1  | 54 | 25 | 32 | 5  | 36 | 23 | 50 | 19 | 44 | 9  | 40 | 15 | 46 | **399** |
| 28 | 31 | 4  | 53 | 24 | 51 | 20 | 35 | 8  | 39 | 14 | 45 | 10 | 41 | 403 |
| 55 | 2  | 29 | 26 | 33 | 6  | 37 | 22 | 49 | 18 | 43 | 12 | 47 | 16 | 395 |
| 30 | 27 | 56 | 3  | 52 | 21 | 34 | 7  | 38 | 13 | 48 | 17 | 42 | 11 | **399** |
| **114** | **114** | **114** | **114** | **114** | **114** | **114** | **114** | **114** | **114** | **114** | **114** | **114** | **114** | |

Fig.12. Near-magic tour on 4 x 14 board with short rows adding to magic constant 114.

There are 1092618 semi-magic tours, 170 quasi-magic tours and 330 near-magic tours with long rows adding to magic constant 399. They are shown in Figure 13 to Figure 15.

| 1  | 46 | 3  | 52 | 11 | 48 | 17 | 50 | 13 | 40 | 27 | 32 | 23 | 36 | **399** |
| 4  | 53 | 10 | 47 | 8  | 51 | 12 | 41 | 28 | 31 | 20 | 35 | 26 | 33 | **399** |
| 45 | 2  | 55 | 6  | 43 | 16 | 49 | 18 | 39 | 14 | 29 | 24 | 37 | 22 | **399** |
| 54 | 5  | 44 | 9  | 56 | 7  | 42 | 15 | 30 | 19 | 38 | 21 | 34 | 25 | **399** |
| 104 | 106 | 112 | **114** | 118 | 122 | 120 | 124 | 110 | 104 | **114** | 112 | 120 | 116 | |

Fig.13. Semi-magic tour on 4 x 14 board with long rows adding to magic constant 399.

| 3  | 52 | 19 | 40 | 1  | 50 | 21 | 42 | 13 | 48 | 23 | 30 | 11 | 46 | **399** |
| 18 | 39 | 2  | 51 | 20 | 41 | 14 | 49 | 22 | 29 | 12 | 47 | 24 | 31 | **399** |
| 53 | 4  | 37 | 16 | 55 | 6  | 35 | 28 | 43 | 8  | 33 | 26 | 45 | 10 | **399** |
| 38 | 17 | 54 | 5  | 36 | 15 | 56 | 7  | 34 | 27 | 44 | 9  | 32 | 25 | **399** |
| 112 | 112 | 112 | 112 | 112 | 112 | 126 | 126 | 112 | 112 | 112 | 112 | 112 | 112 | |

Fig.14. Quasi-magic tour on 4 x 14 board with long rows adding to magic constant 399.

| 1  | 54 | 25 | 32 | 5  | 52 | 23 | 48 | 9  | 38 | 11 | 40 | 17 | 44 | **399** |
| 28 | 31 | 4  | 53 | 24 | 33 | 6  | 37 | 20 | 47 | 18 | 43 | 14 | 41 | **399** |
| 55 | 2  | 29 | 26 | 51 | 22 | 35 | 8  | 49 | 10 | 39 | 12 | 45 | 16 | **399** |
| 30 | 27 | 56 | 3  | 34 | 7  | 50 | 21 | 36 | 19 | 46 | 15 | 42 | 13 | **399** |
| **114** | **114** | **114** | **114** | **114** | **114** | **114** | **114** | **114** | **114** | **114** | 110 | 118 | **114** | |

Fig.15. Near-magic tour on 4 x 14 board with long rows adding to magic constant 399.

**4 x 15 board**: Tours on 4 x 15 board haven't been investigated earlier. The total number of tours hasn't been counted and the author conjectures that it will be around 2.6 billion. There are 8221 semi-magic tours with short rows magic as shown in Figure 16.

| 1  | 58 | 3  | 34 | 5  | 54 | 9  | 36 | 23 | 40 | 11 | 46 | 17 | 42 | 13 | 392 |
| 30 | 33 | 28 | 57 | 8  | 35 | 6  | 39 | 10 | 37 | 16 | 41 | 12 | 47 | 18 | 417 |
| 59 | 2  | 31 | 4  | 53 | 26 | 55 | 22 | 51 | 24 | 45 | 20 | 49 | 14 | 43 | 498 |
| 32 | 29 | 60 | 27 | 56 | 7  | 52 | 25 | 38 | 21 | 50 | 15 | 44 | 19 | 48 | 523 |
| 122 | 122 | 122 | 122 | 122 | 122 | 122 | 122 | 122 | 122 | 122 | 122 | 122 | 122 | 122 | |

Fig.16. Semi-magic tour on 4 x 15 board with short rows adding to magic constant 122.

Figure 17 has 13 short rows adding to magic constant 122 and the long rows are adding to only two values which are consecutive numbers. There are 22 such tours.

| 1  | 56 | 7  | 58 | 9  | 52 | 27 | 34 | 13 | 44 | 25 | 46 | 21 | 42 | 23 | 458 |
|----|----|----|----|----|----|----|----|----|----|----|----|----|----|----|-----|
| 6  | 59 | 4  | 53 | 28 | 31 | 12 | 49 | 26 | 47 | 16 | 43 | 24 | 39 | 20 | 457 |
| 55 | 2  | 57 | 8  | 51 | 10 | 33 | 30 | 35 | 14 | 45 | 18 | 37 | 22 | 41 | 458 |
| 60 | 5  | 54 | 3  | 32 | 29 | 50 | 11 | 48 | 17 | 36 | 15 | 40 | 19 | 38 | 457 |
| **122** | **122** | **122** | **122** | 120 | **122** | **122** | 124 | **122** | **122** | **122** | **122** | **122** | **122** | **122** | |

Fig.17. The closest approach to quasi-magic tour on 4x15 board.

Figure 18 is the rarest (along with its reverse) of all tours. Here, all the short rows are magic and the long rows are adding to four consecutive numbers, that is, from 456 to 459. There are no quasi-magic tours with short rows magic.

| 3  | 60 | 27 | 34 | 5  | 38 | 25 | 54 | 21 | 52 | 17 | 42 | 13 | 50 | 15 | 456 |
|----|----|----|----|----|----|----|----|----|----|----|----|----|----|----|-----|
| 28 | 31 | 4  | 57 | 26 | 55 | 22 | 37 | 8  | 43 | 20 | 51 | 16 | 47 | 12 | 457 |
| 59 | 2  | 33 | 30 | 35 | 6  | 39 | 24 | 53 | 18 | 41 | 10 | 45 | 14 | 49 | 458 |
| 32 | 29 | 58 | 1  | 56 | 23 | 36 | 7  | 40 | 9  | 44 | 19 | 48 | 11 | 46 | 459 |
| **122** | **122** | **122** | **122** | **122** | **122** | **122** | **122** | **122** | **122** | **122** | **122** | **122** | **122** | **122** | |

Fig.18. The rarest semi-magic tour with short rows magic on 4x15 board.

**4 x 16 board**: Since the number of tours increases rapidly with board size, the author conjectures that it will be over 12 billion. There are 19212 semi-magic tours, 488 quasi-magic tours and 1012 near-magic tours with short rows magic. Such tours are shown in Figure 19 to Figure 21.

| 1  | 62 | 3  | 36 | 29 | 38 | 19 | 58 | 25 | 44 | 17 | 56 | 23 | 52 | 15 | 54 | 532 |
|----|----|----|----|----|----|----|----|----|----|----|----|----|----|----|----|-----|
| 32 | 35 | 30 | 61 | 4  | 59 | 26 | 45 | 18 | 57 | 24 | 43 | 16 | 55 | 12 | 51 | 568 |
| 63 | 2  | 33 | 28 | 37 | 6  | 39 | 20 | 47 | 8  | 41 | 22 | 49 | 10 | 53 | 14 | 472 |
| 34 | 31 | 64 | 5  | 60 | 27 | 46 | 7  | 40 | 21 | 48 | 9  | 42 | 13 | 50 | 11 | 508 |
| **130** | **130** | **130** | **130** | **130** | **130** | **130** | **130** | **130** | **130** | **130** | **130** | **130** | **130** | **130** | **130** | |

Fig.19. Semi-magic tour on 4 x 16 board with short rows adding to magic constant 130.

| 1  | 62 | 31 | 36 | 27 | 38 | 7  | 42 | 25 | 56 | 21 | 50 | 11 | 46 | 17 | 52 | 522 |
|----|----|----|----|----|----|----|----|----|----|----|----|----|----|----|----|-----|
| 32 | 35 | 2  | 61 | 6  | 59 | 26 | 57 | 22 | 41 | 10 | 45 | 16 | 51 | 12 | 47 | 522 |
| 63 | 30 | 33 | 4  | 37 | 28 | 39 | 8  | 43 | 24 | 55 | 20 | 49 | 14 | 53 | 18 | 518 |
| 34 | 3  | 64 | 29 | 60 | 5  | 58 | 23 | 40 | 9  | 44 | 15 | 54 | 19 | 48 | 13 | 518 |
| **130** | **130** | **130** | **130** | **130** | **130** | **130** | **130** | **130** | **130** | **130** | **130** | **130** | **130** | **130** | **130** | |

Fig.20. Quasi-magic tour on 4 x 16 board with short rows adding to magic constant 130.

| 1  | 62 | 3  | 36 | 27 | 60 | 7  | 56 | 25 | 42 | 21 | 54 | 13 | 46 | 17 | 50 | **520** |
|----|----|----|----|----|----|----|----|----|----|----|----|----|----|----|----|---------|
| 32 | 35 | 30 | 61 | 6  | 37 | 26 | 41 | 8  | 55 | 12 | 43 | 20 | 49 | 14 | 47 | 516 |
| 63 | 2  | 33 | 4  | 59 | 28 | 39 | 24 | 57 | 10 | 53 | 22 | 45 | 16 | 51 | 18 | 524 |
| 34 | 31 | 64 | 29 | 38 | 5  | 58 | 9  | 40 | 23 | 44 | 11 | 52 | 19 | 48 | 15 | **520** |
| **130** | **130** | **130** | **130** | **130** | **130** | **130** | **130** | **130** | **130** | **130** | **130** | **130** | **130** | **130** | **130** | |

Fig.21. Near-magic tour on 4 x 16 board with short rows adding to magic constant 130.

The number of semi-magic tours with long rows magic has not been enumerated. There are 710 quasi-magic tours and 1304 near-magic tours with long rows magic. They are shown in Figure 22 to Figure 24.

| 1 | 52 | 3 | 62 | 11 | 54 | 15 | 58 | 19 | 46 | 23 | 40 | 31 | 36 | 27 | 42 | **520** |
|---|---|---|---|---|---|---|---|---|---|---|---|---|---|---|---|---|
| 4 | 63 | 12 | 53 | 14 | 59 | 8 | 55 | 16 | 39 | 32 | 37 | 22 | 41 | 30 | 35 | **520** |
| 51 | 2 | 61 | 6 | 49 | 10 | 57 | 18 | 47 | 20 | 45 | 24 | 33 | 28 | 43 | 26 | **520** |
| 64 | 5 | 50 | 13 | 60 | 7 | 48 | 9 | 56 | 17 | 38 | 21 | 44 | 25 | 34 | 29 | **520** |
| 120 | 122 | 126 | 134 | 134 | **130** | 128 | 140 | 138 | 122 | 138 | 122 | **130** | **130** | 134 | 132 | |

Fig.22. Semi-magic tour on 4 x 16 board with long rows adding to magic constant 520.

| 5 | 52 | 11 | 62 | 3 | 50 | 9 | 60 | 17 | 40 | 27 | 46 | 31 | 38 | 25 | 44 | **520** |
|---|---|---|---|---|---|---|---|---|---|---|---|---|---|---|---|---|
| 12 | 55 | 4 | 51 | 10 | 61 | 2 | 49 | 28 | 47 | 32 | 39 | 26 | 45 | 22 | 37 | **520** |
| 53 | 6 | 57 | 14 | 63 | 8 | 59 | 16 | 33 | 18 | 41 | 30 | 35 | 20 | 43 | 24 | **520** |
| 56 | 13 | 54 | 7 | 58 | 15 | 64 | 1 | 48 | 29 | 34 | 19 | 42 | 23 | 36 | 21 | **520** |
| 126 | 126 | 126 | 134 | 134 | 134 | 134 | 126 | 126 | 134 | 134 | 134 | 134 | 126 | 126 | 126 | |

Fig.23. Quasi-magic tour on 4 x 16 board with long rows adding to magic constant 520.

| 3 | 60 | 7 | 64 | 21 | 56 | 11 | 42 | 19 | 52 | 25 | 48 | 15 | 34 | 27 | 36 | **520** |
|---|---|---|---|---|---|---|---|---|---|---|---|---|---|---|---|---|
| 6 | 63 | 4 | 57 | 10 | 43 | 20 | 53 | 24 | 47 | 16 | 51 | 26 | 37 | 30 | 33 | **520** |
| 59 | 2 | 61 | 8 | 55 | 22 | 45 | 12 | 41 | 18 | 49 | 14 | 39 | 32 | 35 | 28 | **520** |
| 62 | 5 | 58 | 1 | 44 | 9 | 54 | 23 | 46 | 13 | 40 | 17 | 50 | 29 | 38 | 31 | **520** |
| **130** | **130** | **130** | **130** | **130** | **130** | **130** | **130** | **130** | **130** | **130** | **130** | **130** | 132 | **130** | 128 | |

Fig.24. Near-magic tour on 4 x 16 board with long rows adding to magic constant 520.

**4 x 17 board**: The total number of tours hasn't been computed. There are 45262 semi-magic tours with short rows magic. Figure 25 shows one such tour. Discerning readers must have spotted that the sum of four long rows are consecutive numbers from 585 to 588. There are 30 such semi-magic tours. Figure 26 shows two interesting tours. They have fifteen short rows adding to magic constant 138 and the remaining two rows are +-2 from the magic constant. Long rows are adding to consecutive numbers, 586 and 587. There are only 51 tours with these properties and no quasi-magic tours with short rows magic.

| 1 | 34 | 67 | 36 | **138** |
|---|---|---|---|---|
| 66 | 37 | 2 | 33 | **138** |
| 31 | 4 | 35 | 68 | **138** |
| 38 | 65 | 32 | 3 | **138** |
| 5 | 30 | 39 | 64 | **138** |
| 42 | 63 | 6 | 27 | **138** |
| 29 | 26 | 43 | 40 | **138** |
| 62 | 41 | 28 | 7 | **138** |
| 25 | 8 | 61 | 44 | **138** |
| 60 | 45 | 24 | 9 | **138** |
| 23 | 12 | 57 | 46 | **138** |
| 56 | 59 | 10 | 13 | **138** |
| 11 | 22 | 47 | 58 | **138** |
| 48 | 55 | 14 | 21 | **138** |
| 17 | 20 | 49 | 52 | **138** |
| 54 | 51 | 18 | 15 | **138** |
| 19 | 16 | 53 | 50 | **138** |
| 587 | 588 | 585 | 586 | |

Fig.25. Semi-magic tour

| 1 | 4 | 65 | 68 | **138** |
|---|---|---|---|---|
| 64 | 67 | 2 | 5 | **138** |
| 3 | 10 | 59 | 66 | **138** |
| 60 | 63 | 6 | 9 | **138** |
| 11 | 8 | 61 | 58 | **138** |
| 62 | 57 | 12 | 7 | **138** |
| 13 | 34 | 55 | 36 | **138** |
| 56 | 37 | 32 | 15 | 140 |
| 33 | 14 | 35 | 54 | 136 |
| 38 | 53 | 16 | 31 | **138** |
| 23 | 30 | 39 | 46 | **138** |
| 52 | 45 | 24 | 17 | **138** |
| 29 | 22 | 47 | 40 | **138** |
| 44 | 51 | 18 | 25 | **138** |
| 21 | 28 | 41 | 48 | **138** |
| 50 | 43 | 26 | 19 | **138** |
| 27 | 20 | 49 | 42 | **138** |
| 587 | 586 | 587 | 586 | |

| 1 | 4 | 65 | 68 | **138** |
|---|---|---|---|---|
| 64 | 67 | 2 | 5 | **138** |
| 3 | 10 | 59 | 66 | **138** |
| 60 | 63 | 6 | 9 | **138** |
| 11 | 8 | 61 | 58 | **138** |
| 62 | 57 | 12 | 7 | **138** |
| 13 | 34 | 55 | 36 | **138** |
| 56 | 37 | 32 | 15 | 140 |
| 33 | 14 | 35 | 54 | 136 |
| 38 | 53 | 16 | 31 | **138** |
| 27 | 30 | 39 | 42 | **138** |
| 52 | 41 | 28 | 17 | **138** |
| 29 | 26 | 43 | 40 | **138** |
| 44 | 51 | 18 | 25 | **138** |
| 21 | 24 | 45 | 48 | **138** |
| 50 | 47 | 22 | 19 | **138** |
| 23 | 20 | 49 | 46 | **138** |
| 587 | 586 | 587 | 586 | |

Fig.26. Tours with two short rows +-2 to magic constant 138.

**4 x 18 board**: The total number of tours hasn't been counted. There are 213280 semi-magic tours, 2624 quasi-magic tours and 3976 near-magic tours with short rows magic. There are 2492 quasi-magic tours and 6322 near-magic tours with long rows magic. The author is not giving any example of these tours because they have eclipsed in the light of long awaited and much sought after MAGIC TOUR. *There are 16 magic tours on 4 x 18 board.* They are shown in Figure 27. All short rows sum to 146 and the long rows to 657.

| 1 | 70 | 33 | 40 | 5 | 42 | 9 | 66 | 29 | 62 | 27 | 58 | 13 | 52 | 25 | 56 | 19 | 50 |
|---|---|---|---|---|---|---|---|---|---|---|---|---|---|---|---|---|---|
| 36 | 39 | 4 | 69 | 32 | 65 | 6 | 43 | 10 | 45 | 14 | 61 | 26 | 57 | 20 | 51 | 24 | 55 |
| 71 | 2 | 37 | 34 | 41 | 8 | 67 | 30 | 63 | 28 | 59 | 12 | 47 | 16 | 53 | 22 | 49 | 18 |
| 38 | 35 | 72 | 3 | 68 | 31 | 64 | 7 | 44 | 11 | 46 | 15 | 60 | 21 | 48 | 17 | 54 | 23 |

| 1 | 70 | 33 | 40 | 5 | 42 | 9 | 66 | 29 | 62 | 13 | 58 | 27 | 56 | 25 | 52 | 19 | 50 |
|---|---|---|---|---|---|---|---|---|---|---|---|---|---|---|---|---|---|
| 36 | 39 | 4 | 69 | 32 | 65 | 6 | 43 | 10 | 59 | 28 | 61 | 16 | 47 | 18 | 49 | 22 | 53 |
| 71 | 2 | 37 | 34 | 41 | 8 | 67 | 30 | 63 | 14 | 45 | 12 | 57 | 26 | 55 | 24 | 51 | 20 |
| 38 | 35 | 72 | 3 | 68 | 31 | 64 | 7 | 44 | 11 | 60 | 15 | 46 | 17 | 48 | 21 | 54 | 23 |

| 1 | 70 | 33 | 40 | 5 | 44 | 31 | 66 | 27 | 50 | 9 | 48 | 21 | 62 | 19 | 60 | 15 | 56 |
|---|---|---|---|---|---|---|---|---|---|---|---|---|---|---|---|---|---|
| 36 | 39 | 4 | 69 | 32 | 67 | 28 | 43 | 8 | 47 | 24 | 51 | 10 | 53 | 12 | 57 | 18 | 59 |
| 71 | 2 | 37 | 34 | 41 | 6 | 45 | 30 | 65 | 26 | 49 | 22 | 63 | 20 | 61 | 16 | 55 | 14 |
| 38 | 35 | 72 | 3 | 68 | 29 | 42 | 7 | 46 | 23 | 64 | 25 | 52 | 11 | 54 | 13 | 58 | 17 |

| 1 | 70 | 33 | 40 | 5 | 44 | 31 | 66 | 27 | 64 | 23 | 48 | 11 | 52 | 17 | 58 | 13 | 54 |
|---|---|---|---|---|---|---|---|---|---|---|---|---|---|---|---|---|---|
| 36 | 39 | 4 | 69 | 32 | 67 | 28 | 43 | 8 | 47 | 10 | 51 | 24 | 57 | 12 | 53 | 18 | 59 |
| 71 | 2 | 37 | 34 | 41 | 6 | 45 | 30 | 65 | 26 | 63 | 22 | 49 | 16 | 61 | 20 | 55 | 14 |
| 38 | 35 | 72 | 3 | 68 | 29 | 42 | 7 | 46 | 9 | 50 | 25 | 62 | 21 | 56 | 15 | 60 | 19 |

| 3 | 72 | 33 | 40 | 5 | 42 | 9 | 66 | 29 | 46 | 27 | 60 | 25 | 50 | 21 | 58 | 17 | 54 |
|---|---|---|---|---|---|---|---|---|---|---|---|---|---|---|---|---|---|
| 34 | 37 | 4 | 69 | 32 | 65 | 6 | 43 | 10 | 61 | 12 | 47 | 22 | 59 | 24 | 55 | 20 | 57 |
| 71 | 2 | 39 | 36 | 41 | 8 | 67 | 30 | 63 | 28 | 45 | 26 | 51 | 14 | 49 | 18 | 53 | 16 |
| 38 | 35 | 70 | 1 | 68 | 31 | 64 | 7 | 44 | 11 | 62 | 13 | 48 | 23 | 52 | 15 | 56 | 19 |

| 3 | 72 | 33 | 40 | 5 | 42 | 9 | 66 | 29 | 62 | 11 | 60 | 25 | 50 | 21 | 58 | 17 | 54 |
|---|---|---|---|---|---|---|---|---|---|---|---|---|---|---|---|---|---|
| 34 | 37 | 4 | 69 | 32 | 65 | 6 | 43 | 10 | 45 | 28 | 47 | 22 | 59 | 24 | 55 | 20 | 57 |
| 71 | 2 | 39 | 36 | 41 | 8 | 67 | 30 | 63 | 12 | 61 | 26 | 51 | 14 | 49 | 18 | 53 | 16 |
| 38 | 35 | 70 | 1 | 68 | 31 | 64 | 7 | 44 | 27 | 46 | 13 | 48 | 23 | 52 | 15 | 56 | 19 |

| 3 | 72 | 33 | 40 | 5 | 44 | 31 | 66 | 27 | 48 | 25 | 62 | 15 | 50 | 13 | 54 | 17 | 52 |
|---|---|---|---|---|---|---|---|---|---|---|---|---|---|---|---|---|---|
| 34 | 37 | 4 | 69 | 32 | 67 | 28 | 43 | 8 | 63 | 10 | 49 | 12 | 59 | 16 | 51 | 20 | 55 |
| 71 | 2 | 39 | 36 | 41 | 6 | 45 | 30 | 65 | 26 | 47 | 24 | 61 | 14 | 57 | 22 | 53 | 18 |
| 38 | 35 | 70 | 1 | 68 | 29 | 42 | 7 | 46 | 9 | 64 | 11 | 58 | 23 | 60 | 19 | 56 | 21 |

| 3  | 72 | 33 | 40 | 5  | 44 | 31 | 66 | 27 | 64 | 9  | 62 | 15 | 50 | 13 | 54 | 17 | 52 |
|----|----|----|----|----|----|----|----|----|----|----|----|----|----|----|----|----|----|
| 34 | 37 | 4  | 69 | 32 | 67 | 28 | 43 | 8  | 47 | 26 | 49 | 12 | 59 | 16 | 51 | 20 | 55 |
| 71 | 2  | 39 | 36 | 41 | 6  | 45 | 30 | 65 | 10 | 63 | 24 | 61 | 14 | 57 | 22 | 53 | 18 |
| 38 | 35 | 70 | 1  | 68 | 29 | 42 | 7  | 46 | 25 | 48 | 11 | 58 | 23 | 60 | 19 | 56 | 21 |

| 17 | 58 | 13 | 54 | 11 | 52 | 25 | 64 | 23 | 46 | 7  | 42 | 29 | 68 | 33 | 40 | 3  | 72 |
|----|----|----|----|----|----|----|----|----|----|----|----|----|----|----|----|----|----|
| 14 | 55 | 16 | 61 | 20 | 63 | 22 | 49 | 26 | 65 | 30 | 45 | 6  | 41 | 4  | 69 | 34 | 37 |
| 59 | 18 | 57 | 12 | 53 | 10 | 51 | 24 | 47 | 8  | 43 | 28 | 67 | 32 | 39 | 36 | 71 | 2  |
| 56 | 15 | 60 | 19 | 62 | 21 | 48 | 9  | 50 | 27 | 66 | 31 | 44 | 5  | 70 | 1  | 38 | 35 |

| 19 | 60 | 15 | 56 | 21 | 62 | 25 | 50 | 9  | 46 | 7  | 42 | 29 | 68 | 33 | 40 | 3  | 72 |
|----|----|----|----|----|----|----|----|----|----|----|----|----|----|----|----|----|----|
| 14 | 55 | 20 | 61 | 16 | 49 | 22 | 63 | 26 | 65 | 30 | 45 | 6  | 41 | 4  | 69 | 34 | 37 |
| 59 | 18 | 53 | 12 | 57 | 24 | 51 | 10 | 47 | 8  | 43 | 28 | 67 | 32 | 39 | 36 | 71 | 2  |
| 54 | 13 | 58 | 17 | 52 | 11 | 48 | 23 | 64 | 27 | 66 | 31 | 44 | 5  | 70 | 1  | 38 | 35 |

| 19 | 56 | 15 | 52 | 23 | 48 | 13 | 46 | 27 | 44 | 7  | 64 | 31 | 68 | 33 | 40 | 1  | 70 |
|----|----|----|----|----|----|----|----|----|----|----|----|----|----|----|----|----|----|
| 16 | 53 | 18 | 49 | 14 | 51 | 26 | 61 | 12 | 63 | 30 | 67 | 8  | 41 | 4  | 69 | 36 | 39 |
| 57 | 20 | 55 | 24 | 59 | 22 | 47 | 28 | 45 | 10 | 43 | 6  | 65 | 32 | 37 | 34 | 71 | 2  |
| 54 | 17 | 58 | 21 | 50 | 25 | 60 | 11 | 62 | 29 | 66 | 9  | 42 | 5  | 72 | 3  | 38 | 35 |

| 19 | 56 | 15 | 52 | 23 | 48 | 13 | 62 | 11 | 44 | 7  | 64 | 31 | 68 | 33 | 40 | 1  | 70 |
|----|----|----|----|----|----|----|----|----|----|----|----|----|----|----|----|----|----|
| 16 | 53 | 18 | 49 | 14 | 51 | 26 | 45 | 28 | 63 | 30 | 67 | 8  | 41 | 4  | 69 | 36 | 39 |
| 57 | 20 | 55 | 24 | 59 | 22 | 47 | 12 | 61 | 10 | 43 | 6  | 65 | 32 | 37 | 34 | 71 | 2  |
| 54 | 17 | 58 | 21 | 50 | 25 | 60 | 27 | 46 | 29 | 66 | 9  | 42 | 5  | 72 | 3  | 38 | 35 |

| 21 | 56 | 19 | 60 | 23 | 58 | 11 | 48 | 25 | 46 | 7  | 42 | 29 | 68 | 33 | 40 | 1  | 70 |
|----|----|----|----|----|----|----|----|----|----|----|----|----|----|----|----|----|----|
| 18 | 53 | 22 | 57 | 14 | 61 | 24 | 63 | 10 | 65 | 30 | 45 | 6  | 41 | 4  | 69 | 36 | 39 |
| 55 | 20 | 51 | 16 | 59 | 12 | 49 | 26 | 47 | 8  | 43 | 28 | 67 | 32 | 37 | 34 | 71 | 2  |
| 52 | 17 | 54 | 13 | 50 | 15 | 62 | 9  | 64 | 27 | 66 | 31 | 44 | 5  | 72 | 3  | 38 | 35 |

| 21 | 56 | 19 | 60 | 23 | 58 | 11 | 64 | 9  | 46 | 7  | 42 | 29 | 68 | 33 | 40 | 1  | 70 |
|----|----|----|----|----|----|----|----|----|----|----|----|----|----|----|----|----|----|
| 18 | 53 | 22 | 57 | 14 | 61 | 24 | 47 | 26 | 65 | 30 | 45 | 6  | 41 | 4  | 69 | 36 | 39 |
| 55 | 20 | 51 | 16 | 59 | 12 | 49 | 10 | 63 | 8  | 43 | 28 | 67 | 32 | 37 | 34 | 71 | 2  |
| 52 | 17 | 54 | 13 | 50 | 15 | 62 | 25 | 48 | 27 | 66 | 31 | 44 | 5  | 72 | 3  | 38 | 35 |

| 23 | 54 | 17 | 48 | 21 | 60 | 15 | 46 | 11 | 44 | 7  | 64 | 31 | 68 | 33 | 40 | 3  | 72 |
|----|----|----|----|----|----|----|----|----|----|----|----|----|----|----|----|----|----|
| 18 | 49 | 22 | 53 | 16 | 47 | 12 | 59 | 28 | 63 | 30 | 67 | 8  | 41 | 4  | 69 | 34 | 37 |
| 55 | 24 | 51 | 20 | 57 | 26 | 61 | 14 | 45 | 10 | 43 | 6  | 65 | 32 | 39 | 36 | 71 | 2  |
| 50 | 19 | 56 | 25 | 52 | 13 | 58 | 27 | 62 | 29 | 66 | 9  | 42 | 5  | 70 | 1  | 38 | 35 |

| 23 | 54 | 21 | 48 | 17 | 46 | 15 | 60 | 11 | 44 | 7  | 64 | 31 | 68 | 33 | 40 | 3  | 72 |
|----|----|----|----|----|----|----|----|----|----|----|----|----|----|----|----|----|----|
| 20 | 51 | 24 | 55 | 26 | 57 | 12 | 45 | 14 | 63 | 30 | 67 | 8  | 41 | 4  | 69 | 34 | 37 |
| 53 | 22 | 49 | 18 | 47 | 16 | 61 | 28 | 59 | 10 | 43 | 6  | 65 | 32 | 39 | 36 | 71 | 2  |
| 50 | 19 | 52 | 25 | 56 | 27 | 58 | 13 | 62 | 29 | 66 | 9  | 42 | 5  | 70 | 1  | 38 | 35 |

Fig.27. **Magic tours on 4 x 18 board.**

**4 x 19 board**: There are 250247 semi-magic tours with short rows magic as shown in Figure 28. The tour in Figure 29 has 17 short rows magic and long rows adding to two different values. There are 16042 such tours. There are no quasi-magic tours with short rows magic.

| 1 | 74 | 3 | 42 | 7 | 72 | 33 | 68 | 31 | 66 | 13 | 50 | 27 | 62 | 25 | 56 | 19 | 60 | 23 | 732 |
|---|---|---|---|---|---|---|---|---|---|---|---|---|---|---|---|---|---|---|---|
| 38 | 41 | 36 | 73 | 34 | 69 | 6 | 45 | 10 | 47 | 28 | 63 | 14 | 51 | 20 | 61 | 24 | 55 | 18 | 733 |
| 75 | 2 | 39 | 4 | 43 | 8 | 71 | 32 | 67 | 30 | 65 | 12 | 49 | 26 | 57 | 16 | 53 | 22 | 59 | 730 |
| 40 | 37 | 76 | 35 | 70 | 5 | 44 | 9 | 46 | 11 | 48 | 29 | 64 | 15 | 52 | 21 | 58 | 17 | 54 | 731 |
| **154** | **154** | **154** | **154** | **154** | **154** | **154** | **154** | **154** | **154** | **154** | **154** | **154** | **154** | **154** | **154** | **154** | **154** | **154** | |

Fig.28. Semi-magic tour on 4 x 19 board with short rows magic.

| 1 | 74 | 5 | 70 | 11 | 68 | 13 | 64 | 35 | 42 | 17 | 56 | 33 | 58 | 27 | 54 | 31 | 48 | 25 | 732 |
|---|---|---|---|---|---|---|---|---|---|---|---|---|---|---|---|---|---|---|---|
| 4 | 71 | 2 | 67 | 8 | 65 | 38 | 41 | 16 | 61 | 34 | 59 | 20 | 55 | 32 | 49 | 26 | 53 | 30 | 731 |
| 73 | 6 | 75 | 10 | 69 | 12 | 63 | 14 | 39 | 36 | 43 | 18 | 57 | 22 | 45 | 28 | 51 | 24 | 47 | 732 |
| 76 | 3 | 72 | 7 | 66 | 9 | 40 | 37 | 62 | 15 | 60 | 21 | 44 | 19 | 50 | 23 | 46 | 29 | 52 | 731 |
| **154** | **154** | **154** | **154** | **154** | **154** | **154** | 156 | 152 | **154** | **154** | **154** | **154** | **154** | **154** | **154** | **154** | **154** | **154** | |

Fig. 29. The closest approach to quasi-magic tour on 4 x 19 board.

**4 x 20 board**: The total number of tours hasn't been counted. There are 587072 semi-magic tours, 12420 quasi-magic tours and 18440 near-magic tours with short rows magic. There are 88 magic tours and three of them are shown in Figure 30.

| 1 | 40 | 79 | 42 | **162** | 1 | 40 | 79 | 42 | **162** | 1 | 40 | 79 | 42 | **162** |
|---|---|---|---|---|---|---|---|---|---|---|---|---|---|---|
| 78 | 43 | 2 | 39 | **162** | 78 | 43 | 2 | 39 | **162** | 78 | 43 | 2 | 39 | **162** |
| 3 | 38 | 41 | 80 | **162** | 3 | 38 | 41 | 80 | **162** | 3 | 38 | 41 | 80 | **162** |
| 44 | 77 | 4 | 37 | **162** | 44 | 77 | 4 | 37 | **162** | 44 | 77 | 4 | 37 | **162** |
| 5 | 34 | 47 | 76 | **162** | 7 | 36 | 45 | 74 | **162** | 35 | 6 | 75 | 46 | **162** |
| 48 | 45 | 36 | 33 | **162** | 76 | 73 | 8 | 5 | **162** | 76 | 45 | 36 | 5 | **162** |
| 35 | 6 | 75 | 46 | **162** | 35 | 6 | 75 | 46 | **162** | 7 | 34 | 47 | 74 | **162** |
| 74 | 49 | 32 | 7 | **162** | 72 | 47 | 34 | 9 | **162** | 48 | 71 | 10 | 33 | **162** |
| 31 | 8 | 73 | 50 | **162** | 33 | 10 | 71 | 48 | **162** | 11 | 8 | 73 | 70 | **162** |
| 72 | 51 | 30 | 9 | **162** | 70 | 49 | 32 | 11 | **162** | 72 | 49 | 32 | 9 | **162** |
| 29 | 10 | 71 | 52 | **162** | 31 | 12 | 69 | 50 | **162** | 31 | 12 | 69 | 50 | **162** |
| 70 | 55 | 26 | 11 | **162** | 66 | 51 | 30 | 15 | **162** | 68 | 51 | 30 | 13 | **162** |
| 25 | 28 | 53 | 56 | **162** | 13 | 16 | 65 | 68 | **162** | 29 | 26 | 55 | 52 | **162** |
| 54 | 69 | 12 | 27 | **162** | 52 | 67 | 14 | 29 | **162** | 54 | 67 | 14 | 27 | **162** |
| 13 | 24 | 57 | 68 | **162** | 17 | 28 | 53 | 64 | **162** | 25 | 28 | 53 | 56 | **162** |
| 62 | 67 | 14 | 19 | **162** | 54 | 59 | 22 | 27 | **162** | 66 | 57 | 24 | 15 | **162** |
| 23 | 18 | 63 | 58 | **162** | 23 | 18 | 63 | 58 | **162** | 23 | 16 | 65 | 58 | **162** |
| 66 | 61 | 20 | 15 | **162** | 60 | 55 | 26 | 21 | **162** | 62 | 59 | 22 | 19 | **162** |
| 17 | 22 | 59 | 64 | **162** | 19 | 24 | 57 | 62 | **162** | 17 | 20 | 61 | 64 | **162** |
| 60 | 65 | 16 | 21 | **162** | 56 | 61 | 20 | 25 | **162** | 60 | 63 | 18 | 21 | **162** |
| **810** | **810** | **810** | **810** | | **810** | **810** | **810** | **810** | | **810** | **810** | **810** | **810** | |

Fig.30. Magic tours of knight on 4 x 20 board.

**4 x 21 board**: The total number of tours hasn't been counted. There are 1378912 semi-magic tours with short rows magic as shown in Figure 31. The tours in Figure 32 have 19 short rows magic and long rows are alternately equal with a difference of one. There are 3294 such tours. There are no quasi-magic tours with short rows magic.

| 1 | 42 | 83 | 44 | **170** |
|---|---|---|---|---|
| 82 | 45 | 2 | 41 | **170** |
| 3 | 40 | 43 | 84 | **170** |
| 46 | 81 | 4 | 39 | **170** |
| 5 | 8 | 77 | 80 | **170** |
| 78 | 47 | 38 | 7 | **170** |
| 9 | 6 | 79 | 76 | **170** |
| 48 | 51 | 34 | 37 | **170** |
| 35 | 10 | 75 | 50 | **170** |
| 52 | 49 | 36 | 33 | **170** |
| 11 | 14 | 71 | 74 | **170** |
| 72 | 53 | 32 | 13 | **170** |
| 15 | 12 | 73 | 70 | **170** |
| 54 | 57 | 28 | 31 | **170** |
| 29 | 16 | 69 | 56 | **170** |
| 58 | 55 | 30 | 27 | **170** |
| 17 | 22 | 63 | 68 | **170** |
| 64 | 59 | 26 | 21 | **170** |
| 23 | 18 | 67 | 62 | **170** |
| 60 | 65 | 20 | 25 | **170** |
| 19 | 24 | 61 | 66 | **170** |
| 781 | 774 | 1011 | 1004 | |

Fig.31. Semi-magic tour

| 1 | 4 | 81 | 84 | **170** |
|---|---|---|---|---|
| 80 | 83 | 2 | 5 | **170** |
| 3 | 6 | 79 | 82 | **170** |
| 44 | 77 | 40 | 7 | 168 |
| 41 | 8 | 43 | 78 | **170** |
| 76 | 45 | 10 | 39 | **170** |
| 9 | 42 | 75 | 46 | 172 |
| 48 | 73 | 38 | 11 | **170** |
| 37 | 12 | 47 | 74 | **170** |
| 72 | 49 | 36 | 13 | **170** |
| 35 | 32 | 53 | 50 | **170** |
| 52 | 71 | 14 | 33 | **170** |
| 31 | 34 | 51 | 54 | **170** |
| 70 | 65 | 20 | 15 | **170** |
| 19 | 30 | 55 | 66 | **170** |
| 64 | 69 | 16 | 21 | **170** |
| 29 | 18 | 67 | 56 | **170** |
| 68 | 63 | 22 | 17 | **170** |
| 25 | 28 | 57 | 60 | **170** |
| 62 | 59 | 26 | 23 | **170** |
| 27 | 24 | 61 | 58 | **170** |
| 893 | 892 | 893 | 892 | |

| 27 | 24 | 61 | 58 | **170** |
|---|---|---|---|---|
| 62 | 59 | 26 | 23 | **170** |
| 25 | 28 | 57 | 60 | **170** |
| 70 | 63 | 22 | 15 | **170** |
| 29 | 16 | 69 | 56 | **170** |
| 64 | 71 | 14 | 21 | **170** |
| 17 | 30 | 55 | 68 | **170** |
| 72 | 65 | 20 | 13 | **170** |
| 31 | 18 | 67 | 54 | **170** |
| 66 | 73 | 12 | 19 | **170** |
| 11 | 32 | 53 | 74 | **170** |
| 52 | 75 | 10 | 33 | **170** |
| 35 | 8 | 51 | 76 | **170** |
| 78 | 43 | 34 | 9 | 164 |
| 7 | 36 | 77 | 50 | **170** |
| 44 | 79 | 42 | 5 | **170** |
| 37 | 6 | 49 | 84 | 176 |
| 80 | 45 | 4 | 41 | **170** |
| 1 | 38 | 83 | 48 | **170** |
| 46 | 81 | 40 | 3 | **170** |
| 39 | 2 | 47 | 82 | **170** |
| 893 | 892 | 893 | 892 | |

Fig.32. 19 short rows are magic and long rows are alternately equal.

**4 x 22 board**: The total number of tours hasn't been counted. There are 3237456 semi-magic tours, 64120 quasi-magic tours and 81252 near-magic tours with short rows magic. There are 464 magic tours and three of them are shown in Figure 33.

| 1 | 44 | 87 | 46 | **178** | 17 | 20 | 69 | 72 | **178** | 27 | 30 | 59 | 62 | **178** |
|---|---|---|---|---|---|---|---|---|---|---|---|---|---|---|
| 86 | 47 | 2 | 43 | **178** | 68 | 71 | 18 | 21 | **178** | 58 | 61 | 28 | 31 | **178** |
| 3 | 42 | 45 | 88 | **178** | 19 | 16 | 73 | 70 | **178** | 29 | 26 | 63 | 60 | **178** |
| 48 | 85 | 4 | 41 | **178** | 74 | 67 | 22 | 15 | **178** | 72 | 57 | 32 | 17 | **178** |
| 5 | 38 | 51 | 84 | **178** | 23 | 28 | 61 | 66 | **178** | 25 | 16 | 73 | 64 | **178** |
| 52 | 49 | 40 | 37 | **178** | 62 | 75 | 14 | 27 | **178** | 56 | 71 | 18 | 33 | **178** |
| 39 | 6 | 83 | 50 | **178** | 29 | 24 | 65 | 60 | **178** | 15 | 24 | 65 | 74 | **178** |
| 82 | 53 | 36 | 7 | **178** | 76 | 63 | 26 | 13 | **178** | 70 | 55 | 34 | 19 | **178** |
| 9 | 34 | 81 | 54 | **178** | 25 | 30 | 59 | 64 | **178** | 23 | 14 | 75 | 66 | **178** |
| 80 | 55 | 8 | 35 | **178** | 58 | 77 | 12 | 31 | **178** | 54 | 69 | 20 | 35 | **178** |
| 33 | 10 | 79 | 56 | **178** | 33 | 10 | 79 | 56 | **178** | 13 | 22 | 67 | 76 | **178** |
| 78 | 57 | 32 | 11 | **178** | 78 | 57 | 32 | 11 | **178** | 68 | 53 | 36 | 21 | **178** |
| 31 | 14 | 75 | 58 | **178** | 9 | 34 | 55 | 80 | **178** | 37 | 12 | 77 | 52 | **178** |
| 74 | 77 | 12 | 15 | **178** | 54 | 81 | 8 | 35 | **178** | 78 | 81 | 8 | 11 | **178** |
| 13 | 30 | 59 | 76 | **178** | 7 | 38 | 51 | 82 | **178** | 9 | 38 | 51 | 80 | **178** |
| 62 | 73 | 16 | 27 | **178** | 50 | 53 | 36 | 39 | **178** | 82 | 79 | 10 | 7 | **178** |
| 29 | 26 | 63 | 60 | **178** | 37 | 6 | 83 | 52 | **178** | 39 | 6 | 83 | 50 | **178** |
| 72 | 61 | 28 | 17 | **178** | 84 | 49 | 40 | 5 | **178** | 84 | 49 | 4 | 41 | **178** |
| 25 | 18 | 71 | 64 | **178** | 41 | 4 | 45 | 88 | **178** | 5 | 40 | 45 | 88 | **178** |
| 70 | 67 | 22 | 19 | **178** | 48 | 85 | 42 | 3 | **178** | 48 | 85 | 42 | 3 | **178** |
| 21 | 24 | 65 | 68 | **178** | 1 | 44 | 87 | 46 | **178** | 1 | 44 | 87 | 46 | **178** |
| 66 | 69 | 20 | 23 | **178** | 86 | 47 | 2 | 43 | **178** | 86 | 47 | 2 | 43 | **178** |
| **979** | **979** | **979** | **979** | | **979** | **979** | **979** | **979** | | **979** | **979** | **979** | **979** | |

Fig.33. Magic tours of knight on 4 x 22 board.

**4 x 24 board**: The total number of tours hasn't been counted. There are 17828024 semi-magic tours, 313584 quasi-magic tours and 388188 near-magic tours with short rows magic. There are 2076 magic tours and three of them are shown in Figure 34.

| 1 | 48 | 95 | 50 | **194** | 19 | 22 | 75 | 78 | **194** | 31 | 28 | 69 | 66 | **194** |
|---|---|---|---|---|---|---|---|---|---|---|---|---|---|---|
| 94 | 51 | 2 | 47 | **194** | 74 | 77 | 20 | 23 | **194** | 70 | 67 | 30 | 27 | **194** |
| 3 | 46 | 49 | 96 | **194** | 21 | 18 | 79 | 76 | **194** | 29 | 32 | 65 | 68 | **194** |
| 52 | 93 | 4 | 45 | **194** | 80 | 73 | 24 | 17 | **194** | 64 | 71 | 26 | 33 | **194** |
| 5 | 8 | 89 | 92 | **194** | 25 | 30 | 67 | 72 | **194** | 35 | 24 | 73 | 62 | **194** |
| 90 | 53 | 44 | 7 | **194** | 68 | 81 | 16 | 29 | **194** | 72 | 63 | 34 | 25 | **194** |
| 9 | 6 | 91 | 88 | **194** | 31 | 26 | 71 | 66 | **194** | 23 | 36 | 61 | 74 | **194** |
| 54 | 87 | 10 | 43 | **194** | 82 | 69 | 28 | 15 | **194** | 80 | 75 | 22 | 17 | **194** |
| 13 | 42 | 55 | 84 | **194** | 27 | 32 | 65 | 70 | **194** | 37 | 18 | 79 | 60 | **194** |
| 86 | 83 | 14 | 11 | **194** | 64 | 83 | 14 | 33 | **194** | 76 | 81 | 16 | 21 | **194** |
| 41 | 12 | 85 | 56 | **194** | 13 | 34 | 63 | 84 | **194** | 19 | 38 | 59 | 78 | **194** |
| 82 | 57 | 40 | 15 | **194** | 62 | 85 | 12 | 35 | **194** | 82 | 77 | 20 | 15 | **194** |
| 39 | 16 | 81 | 58 | **194** | 11 | 36 | 61 | 86 | **194** | 39 | 14 | 83 | 58 | **194** |
| 80 | 59 | 38 | 17 | **194** | 60 | 87 | 10 | 37 | **194** | 86 | 57 | 40 | 11 | **194** |
| 37 | 34 | 63 | 60 | **194** | 9 | 40 | 57 | 88 | **194** | 13 | 10 | 87 | 84 | **194** |
| 62 | 79 | 18 | 35 | **194** | 56 | 59 | 38 | 41 | **194** | 56 | 85 | 12 | 41 | **194** |
| 33 | 36 | 61 | 64 | **194** | 39 | 8 | 89 | 58 | **194** | 9 | 42 | 55 | 88 | **194** |
| 78 | 65 | 32 | 19 | **194** | 90 | 55 | 42 | 7 | **194** | 54 | 91 | 6 | 43 | **194** |
| 31 | 20 | 77 | 66 | **194** | 43 | 6 | 91 | 54 | **194** | 5 | 8 | 89 | 92 | **194** |
| 76 | 71 | 26 | 21 | **194** | 52 | 93 | 44 | 5 | **194** | 90 | 53 | 44 | 7 | **194** |
| 25 | 30 | 67 | 72 | **194** | 45 | 4 | 53 | 92 | **194** | 45 | 4 | 93 | 52 | **194** |
| 70 | 75 | 22 | 27 | **194** | 94 | 51 | 48 | 1 | **194** | 50 | 95 | 48 | 1 | **194** |
| 29 | 24 | 73 | 68 | **194** | 3 | 46 | 95 | 50 | **194** | 3 | 46 | 51 | 94 | **194** |
| 74 | 69 | 28 | 23 | **194** | 96 | 49 | 2 | 47 | **194** | 96 | 49 | 2 | 47 | **194** |
| **1164** | **1164** | **1164** | **1164** | | **1164** | **1164** | **1164** | **1164** | | **1164** | **1164** | **1164** | **1164** | |

Fig.34. Magic tours of knight on 4 x 24 board.

**4 x 26 board**: The total number of tours hasn't been counted. There are 98203312 semi-magic tours, 1675648 quasi-magic tours and 1800728 near-magic tours with short rows magic. There are 9904 magic tours and three of them are shown in Figure 35.

| 1 | 52 | 103 | 54 | **210** | 19 | 26 | 79 | 86 | **210** | 35 | 32 | 73 | 70 | **210** |
|---|---|---|---|---|---|---|---|---|---|---|---|---|---|---|
| 102 | 55 | 2 | 51 | **210** | 80 | 87 | 18 | 25 | **210** | 74 | 71 | 34 | 31 | **210** |
| 3 | 50 | 53 | 104 | **210** | 27 | 20 | 85 | 78 | **210** | 33 | 36 | 69 | 72 | **210** |
| 56 | 101 | 4 | 49 | **210** | 88 | 81 | 24 | 17 | **210** | 78 | 75 | 30 | 27 | **210** |
| 5 | 8 | 97 | 100 | **210** | 21 | 28 | 77 | 84 | **210** | 37 | 28 | 77 | 68 | **210** |
| 98 | 57 | 48 | 7 | **210** | 82 | 89 | 16 | 23 | **210** | 76 | 79 | 26 | 29 | **210** |
| 9 | 6 | 99 | 96 | **210** | 29 | 22 | 83 | 76 | **210** | 25 | 38 | 67 | 80 | **210** |
| 58 | 95 | 10 | 47 | **210** | 90 | 75 | 30 | 15 | **210** | 66 | 81 | 24 | 39 | **210** |
| 41 | 46 | 59 | 64 | **210** | 31 | 34 | 71 | 74 | **210** | 21 | 40 | 65 | 84 | **210** |
| 94 | 63 | 42 | 11 | **210** | 72 | 91 | 14 | 33 | **210** | 82 | 85 | 20 | 23 | **210** |
| 45 | 40 | 65 | 60 | **210** | 35 | 32 | 73 | 70 | **210** | 41 | 22 | 83 | 64 | **210** |
| 62 | 93 | 12 | 43 | **210** | 92 | 69 | 36 | 13 | **210** | 86 | 63 | 42 | 19 | **210** |
| 39 | 44 | 61 | 66 | **210** | 37 | 12 | 93 | 68 | **210** | 13 | 18 | 87 | 92 | **210** |
| 92 | 67 | 38 | 13 | **210** | 64 | 67 | 38 | 41 | **210** | 62 | 91 | 14 | 43 | **210** |
| 37 | 14 | 91 | 68 | **210** | 11 | 40 | 65 | 94 | **210** | 17 | 12 | 93 | 88 | **210** |
| 90 | 69 | 36 | 15 | **210** | 66 | 63 | 42 | 39 | **210** | 90 | 61 | 44 | 15 | **210** |
| 17 | 34 | 89 | 70 | **210** | 43 | 10 | 95 | 62 | **210** | 11 | 16 | 89 | 94 | **210** |
| 88 | 71 | 16 | 35 | **210** | 96 | 61 | 8 | 45 | **210** | 60 | 95 | 10 | 45 | **210** |
| 33 | 18 | 87 | 72 | **210** | 9 | 44 | 59 | 98 | **210** | 7 | 46 | 59 | 98 | **210** |
| 86 | 73 | 32 | 19 | **210** | 60 | 97 | 46 | 7 | **210** | 96 | 99 | 6 | 9 | **210** |
| 29 | 20 | 85 | 76 | **210** | 1 | 52 | 99 | 58 | **210** | 47 | 8 | 97 | 58 | **210** |
| 74 | 77 | 28 | 31 | **210** | 100 | 57 | 6 | 47 | **210** | 100 | 57 | 48 | 5 | **210** |
| 21 | 30 | 75 | 84 | **210** | 51 | 2 | 53 | 104 | **210** | 49 | 4 | 53 | 104 | **210** |
| 78 | 81 | 24 | 27 | **210** | 56 | 101 | 48 | 5 | **210** | 56 | 101 | 50 | 3 | **210** |
| 25 | 22 | 83 | 80 | **210** | 3 | 50 | 103 | 54 | **210** | 1 | 52 | 103 | 54 | **210** |
| 82 | 79 | 26 | 23 | **210** | 102 | 55 | 4 | 49 | **210** | 102 | 55 | 2 | 51 | **210** |
| **1365** | **1365** | **1365** | **1365** | | **1365** | **1365** | **1365** | **1365** | | **1365** | **1365** | **1365** | **1365** | |

Fig.35. Magic tours of knight on 4 x 26 board.

**4 x 28 board**: The total number of tours hasn't been counted. There are 540834326 semi-magic tours, 8536072 quasi-magic tours and 8769964 near-magic tours with short rows magic. There are 47456 magic tours and three of them are shown in Figure 36.

| 1 | 56 | 111 | 58 | **226** | 29 | 26 | 87 | 84 | **226** | 35 | 32 | 81 | 78 | **226** |
|---|---|---|---|---|---|---|---|---|---|---|---|---|---|---|
| 110 | 59 | 2 | 55 | **226** | 88 | 85 | 28 | 25 | **226** | 82 | 79 | 34 | 31 | **226** |
| 3 | 54 | 57 | 112 | **226** | 27 | 30 | 83 | 86 | **226** | 33 | 36 | 77 | 80 | **226** |
| 60 | 109 | 4 | 53 | **226** | 82 | 89 | 24 | 31 | **226** | 76 | 83 | 30 | 37 | **226** |
| 5 | 50 | 63 | 108 | **226** | 23 | 32 | 81 | 90 | **226** | 39 | 28 | 75 | 84 | **226** |
| 64 | 61 | 52 | 49 | **226** | 80 | 91 | 34 | 21 | **226** | 86 | 73 | 38 | 29 | **226** |
| 51 | 6 | 107 | 62 | **226** | 33 | 22 | 93 | 78 | **226** | 27 | 40 | 85 | 74 | **226** |
| 106 | 65 | 48 | 7 | **226** | 92 | 79 | 20 | 35 | **226** | 72 | 87 | 42 | 25 | **226** |
| 47 | 8 | 67 | 104 | **226** | 19 | 36 | 77 | 94 | **226** | 41 | 26 | 71 | 88 | **226** |
| 66 | 105 | 46 | 9 | **226** | 76 | 95 | 38 | 17 | **226** | 90 | 69 | 24 | 43 | **226** |
| 45 | 10 | 103 | 68 | **226** | 37 | 18 | 75 | 96 | **226** | 45 | 22 | 89 | 70 | **226** |
| 102 | 69 | 44 | 11 | **226** | 74 | 97 | 16 | 39 | **226** | 68 | 91 | 44 | 23 | **226** |
| 43 | 12 | 101 | 70 | **226** | 15 | 12 | 101 | 98 | **226** | 21 | 46 | 67 | 92 | **226** |
| 100 | 71 | 14 | 41 | **226** | 100 | 73 | 40 | 13 | **226** | 66 | 95 | 18 | 47 | **226** |
| 13 | 42 | 73 | 98 | **226** | 11 | 14 | 99 | 102 | **226** | 17 | 20 | 93 | 96 | **226** |
| 72 | 99 | 40 | 15 | **226** | 72 | 103 | 10 | 41 | **226** | 94 | 65 | 48 | 19 | **226** |
| 39 | 16 | 97 | 74 | **226** | 9 | 44 | 69 | 104 | **226** | 13 | 16 | 97 | 100 | **226** |
| 96 | 93 | 20 | 17 | **226** | 68 | 71 | 42 | 45 | **226** | 64 | 99 | 14 | 49 | **226** |
| 19 | 38 | 75 | 94 | **226** | 43 | 8 | 105 | 70 | **226** | 15 | 12 | 101 | 98 | **226** |
| 92 | 95 | 18 | 21 | **226** | 106 | 67 | 46 | 7 | **226** | 102 | 63 | 50 | 11 | **226** |
| 37 | 24 | 89 | 76 | **226** | 47 | 50 | 63 | 66 | **226** | 51 | 10 | 103 | 62 | **226** |
| 88 | 91 | 22 | 25 | **226** | 64 | 107 | 6 | 49 | **226** | 104 | 107 | 6 | 9 | **226** |
| 23 | 36 | 77 | 90 | **226** | 51 | 48 | 65 | 62 | **226** | 7 | 52 | 61 | 106 | **226** |
| 78 | 87 | 26 | 35 | **226** | 108 | 61 | 52 | 5 | **226** | 108 | 105 | 8 | 5 | **226** |
| 27 | 34 | 79 | 86 | **226** | 53 | 4 | 59 | 110 | **226** | 53 | 4 | 109 | 60 | **226** |
| 80 | 83 | 30 | 33 | **226** | 60 | 109 | 56 | 1 | **226** | 58 | 111 | 56 | 1 | **226** |
| 31 | 28 | 85 | 82 | **226** | 3 | 54 | 111 | 58 | **226** | 3 | 54 | 59 | 110 | **226** |
| 84 | 81 | 32 | 29 | **226** | 112 | 57 | 2 | 55 | **226** | 112 | 57 | 2 | 55 | **226** |
| **1582** | **1582** | **1582** | **1582** | | **1582** | **1582** | **1582** | **1582** | | **1582** | **1582** | **1582** | **1582** | |

Fig.36. Magic tours of knight on 4 x 28 board.

## B. Tours of knight on 6 by n boards

**Introduction**: Even among rectangular boards, tour of knight on 6 by n board has got much less attention vis-à-vis those on 4 by n board. Jelliss [21] has mostly looked into symmetrical tours on 6 by n board. Construction of simple knight's tours on 6 by n (n > 3) boards is relatively easy and the number of closed knight's tours on 6 x 6 board has been counted. But, what about total number of, (both closed and open), knight tours on 6 x 6 board? Can we have magic knight's tours on 6 by n boards? If yes, what is the smallest board size and how many magic tours are there? Magic tours are not possible on a board with odd-even side but can we have quasi-magic tours on such boards? The author plans to look into these questions.

**Description**: A 'magic knight tour' on a rectangular board 6 x n (n even) has all the 6 long (n-cell) rows adding to magic constant $n(6n+1)/2$ and the n short (6-cell) rows adding to magic constant $3(6n+1)$. As is the case with 4 by n board, there can't be magic knight tour on 6 x n board if n is odd.

**6 x 5 board**: Knight tour is not possible on 6 x 2 and 6 x 3 board and tours on 6 x 4 board have already been studied in sufficient detail. Jelliss [21] has only mentioned that there are three closed tours, two of them having axial symmetry, on 6 x 5 board, found by Euler in 1759, Warnsdorf in 1858 and Haldeman in 1864. The author has found that there are 9386 arithmetically distinct tours on 6 x 5 board. Two examples of the semi-magic tours along long (6-cell) rows are shown in Figure 37.

| 1 | 26 | 17 | 22 | 3 | 24 | **93** |
|---|---|---|---|---|---|---|
| 16 | *11* | 2 | 25 | 18 | *21* | **93** |
| 27 | 8 | 19 | 12 | 23 | 4 | **93** |
| *10* | 15 | 6 | 29 | *20* | 13 | **93** |
| 7 | 28 | 9 | 14 | 5 | 30 | **93** |
| 61 | 88 | 53 | 102 | 69 | 92 | |

| 1 | 26 | 17 | 22 | 3 | 24 | **93** |
|---|---|---|---|---|---|---|
| 16 | *21* | 2 | 25 | 18 | *11* | **93** |
| 27 | 8 | 19 | 12 | 23 | 4 | **93** |
| *20* | 15 | 6 | 29 | *10* | 13 | **93** |
| 7 | 28 | 9 | 14 | 5 | 30 | **93** |
| 71 | 98 | 53 | 102 | 59 | 82 | |

Fig.37. Semi-magic knight tours on 6 x 5 board.

Discerning readers must have observed that the two tours are almost identical; only the digits 10, 11 and 20, 21 have swapped places. They have been highlighted in underlined italics. Such semi-magic tours are rare and have aesthetic appeal. Since the reverse of magic/semi-magic tours are also magic/semi-magic, therefore there are four arithmetically different semi-magic tours on 6 x 5 board. There are no quasi-magic tours or near-magic tours with long (6-cell) rows magic. Figure 38 shows two interesting tours. Although they have non-magic lines but the identical sum of rows has an aesthetic appeal. One has four long (6-cell) rows adding up to same value and three short (4-cell) rows also adding up to the same value. The other tour has two long (6-cell) rows adding up to same value, two odd-sum columns and three even-sum columns also adding up to the same value.

| 21 | 8 | 3 | 14 | 19 | 30 | 95 |
|---|---|---|---|---|---|---|
| 2 | 25 | 20 | 29 | 4 | 15 | 95 |
| 9 | 22 | 7 | 16 | 13 | 18 | 85 |
| 26 | 1 | 24 | 11 | 28 | 5 | 95 |
| 23 | 10 | 27 | 6 | 17 | 12 | 95 |
| 81 | 66 | 81 | 76 | 81 | 80 | |

| 27 | 8 | 5 | 18 | 25 | 14 | 97 |
|---|---|---|---|---|---|---|
| 4 | 19 | 26 | 15 | 6 | 17 | 87 |
| 9 | 28 | 7 | 22 | 13 | 24 | 103 |
| 20 | 3 | 30 | 11 | 16 | 1 | 81 |
| 29 | 10 | 21 | 2 | 23 | 12 | 97 |
| 89 | 68 | 89 | 68 | 83 | 68 | |

Fig.38. Knight tour on 6 x 5 board.

**6 x 6 board**: Being a square board, it has received more attention than other 6 x n boards. A closed tour was given by the famous mathematician Leonhard Euler in the year 1759. Jelliss [22] only considered closed tours and has given 9862 distinct 'tour diagrams' computed by Duby in 1964. If we arrange these 'tour diagrams' in numerical form, there are 88758 arithmetically distinct closed tours as per Frenicle's system of classification. Open tours haven't been enumerated earlier and the author has calculated 740982 open tours. So there are 829740 arithmetically distinct tours on 6 x 6 board. Magic tours are always fascinating and Beverley discovered the first one on 8 x 8 board in 1848. So it is natural to look for them on 6 x 6 board too. Before proceeding further, it is better to recollect the following theorem by Jelliss [23]: "A magic knight's tour is impossible on a board with singly-even sides". So magic knight's tours are impossible on boards of size 6 x 6, 6 x 10, 6 x 14 etc. The best we can have is semi-magic tours. Kraitchik discovered the first semi-magic tour in 1927 and Kumar [24] enumerated all 88 arithmetically distinct semi-magic tours in 2002. Figure 39 shows two semi-magic tours. There are no quasi-magic tours or near-magic tours on 6 x 6 board.

| 11 | 36 | 31 | 18 | 13 | 2 | **111** |
|----|----|----|----|----|----|----|
| 32 | 19 | 12 | 1 | 30 | 17 | **111** |
| 35 | 10 | 33 | 16 | 3 | 14 | **111** |
| 20 | 7 | 4 | 27 | 24 | 29 | **111** |
| 5 | 34 | 9 | 22 | 15 | 26 | **111** |
| 8 | 21 | 6 | 25 | 28 | 23 | **111** |
| **111** | 127 | 95 | 109 | 113 | **111** | |

| 17 | 4 | 19 | 36 | 33 | 2 | **111** |
|----|----|----|----|----|----|----|
| 20 | 27 | 16 | 3 | 10 | 35 | **111** |
| 5 | 18 | 21 | 34 | 1 | 32 | **111** |
| 26 | 15 | 28 | 9 | 22 | 11 | **111** |
| 29 | 6 | 13 | 24 | 31 | 8 | **111** |
| 14 | 25 | 30 | 7 | 12 | 23 | **111** |
| **111** | 95 | 127 | 113 | 109 | **111** | |

Fig.39. Semi-magic knight tours on 6 x 6 board.

**6 x 7 board**: Jelliss [25] has given 1067638 distinct closed tour diagrams and the number of open tours is unknown. There are 288 semi-magic tours, all of them open, and two such tours are shown in Figure 40.

| 23 | 32 | 25 | 18 | 21 | 10 | **129** |
|----|----|----|----|----|----|----|
| 26 | 19 | 22 | 11 | 34 | 17 | **129** |
| 31 | 24 | 33 | 20 | 9 | 12 | **129** |
| 4 | 27 | 8 | 39 | 16 | 35 | **129** |
| 7 | 30 | 5 | 36 | 13 | 38 | **129** |
| 42 | 3 | 28 | 15 | 40 | 1 | **129** |
| 29 | 6 | 41 | 2 | 37 | 14 | **129** |
| 162 | 141 | 162 | 141 | 170 | 127 | |

| 23 | 40 | 25 | 18 | 21 | 2 | **129** |
|----|----|----|----|----|----|----|
| 42 | 17 | 22 | 3 | 26 | 19 | **129** |
| 39 | 24 | 41 | 20 | 1 | 4 | **129** |
| 16 | 33 | 38 | 5 | 10 | 27 | **129** |
| 37 | 30 | 35 | 8 | 13 | 6 | **129** |
| 34 | 15 | 32 | 11 | 28 | 9 | **129** |
| 31 | 36 | 29 | 14 | 7 | 12 | **129** |
| 222 | 195 | 222 | 79 | 106 | 79 | |

Fig.40. Semi-magic knight tours on 6 x 7 board.

**6 x 8 board**: Jelliss [25] has given 55488142 distinct closed tour diagrams and the number of open tours is unknown. The author has enumerated 1880 arithmetically distinct semi-magic tours (1572 open and 308 closed tours) with short (6-cell) rows magic and two such tours, one open and the other closed, are shown in Figure 41. There are 42 near-magic (40 open and 2 closed) tours with short rows magic and two examples, one open and the other closed, are shown in Figure 42. The open tour has four long (8-cell) rows summing to magic constant and the two remaining long rows summing to 208 and 184, that is, they are +- 12 from the magic constant 196. This is the only short rows semi-magic tour (along with its reverse) which has four long magic rows.

| 1  | 42 | 11 | 46 | 3  | 44 | **147** |
|----|----|----|----|----|----|---------|
| 10 | 39 | 2  | 43 | 6  | 47 | **147** |
| 37 | 12 | 41 | 8  | 45 | 4  | **147** |
| 40 | 9  | 38 | 5  | 48 | 7  | **147** |
| 21 | 36 | 13 | 28 | 17 | 32 | **147** |
| 24 | 27 | 22 | 31 | 14 | 29 | **147** |
| 35 | 20 | 25 | 16 | 33 | 18 | **147** |
| 26 | 23 | 34 | 19 | 30 | 15 | **147** |
| 194| 208| 186| **196**| **196**| **196** |

| 19 | 36 | 1  | 38 | 21 | 32 | **147** |
|----|----|----|----|----|----|---------|
| 48 | 5  | 20 | 33 | 2  | 39 | **147** |
| 35 | 18 | 37 | 4  | 31 | 22 | **147** |
| 6  | 47 | 34 | 17 | 40 | 3  | **147** |
| 27 | 16 | 41 | 10 | 23 | 30 | **147** |
| 46 | 7  | 28 | 13 | 42 | 11 | **147** |
| 15 | 26 | 9  | 44 | 29 | 24 | **147** |
| 8  | 45 | 14 | 25 | 12 | 43 | **147** |
| 204| 200| 184| 184| 200| 204 |

Fig.41. Semi-magic knight tours with short rows magic on 6 x 8 board.

| 1  | 42 | 11 | 46 | 3  | 44 | **147** |
|----|----|----|----|----|----|---------|
| 10 | 39 | 2  | 43 | 6  | 47 | **147** |
| 41 | 12 | 37 | 8  | 45 | 4  | **147** |
| 38 | 9  | 40 | 5  | 48 | 7  | **147** |
| 21 | 36 | 13 | 28 | 17 | 32 | **147** |
| 24 | 27 | 22 | 31 | 14 | 29 | **147** |
| 35 | 20 | 25 | 16 | 33 | 18 | **147** |
| 26 | 23 | 34 | 19 | 30 | 15 | **147** |
| **196**| 208| 184| **196**| **196**| **196** |

| 31 | 16 | 45 | 12 | 29 | 14 | **147** |
|----|----|----|----|----|----|---------|
| 46 | 1  | 30 | 15 | 44 | 11 | **147** |
| 17 | 32 | 9  | 48 | 13 | 28 | **147** |
| 2  | 47 | 18 | 27 | 10 | 43 | **147** |
| 33 | 26 | 39 | 8  | 19 | 22 | **147** |
| 38 | 3  | 36 | 21 | 42 | 7  | **147** |
| 25 | 34 | 5  | 40 | 23 | 20 | **147** |
| 4  | 37 | 24 | 35 | 6  | 41 | **147** |
| **196**| **196**| 206| 206| 186| 186 |

Fig.42. Near-magic knight tours with short rows magic on 6 x 8 board.

There are 28 quasi-magic (24 open and 4 closed) tours with short rows magic and two examples are shown in Figure 43.

| 1  | 44 | 9  | 40 | 5  | 48 | **147** |
|----|----|----|----|----|----|---------|
| 8  | 39 | 6  | 43 | 10 | 41 | **147** |
| 45 | 2  | 37 | 12 | 47 | 4  | **147** |
| 38 | 7  | 46 | 3  | 42 | 11 | **147** |
| 17 | 32 | 21 | 36 | 13 | 28 | **147** |
| 20 | 35 | 18 | 27 | 22 | 25 | **147** |
| 31 | 16 | 33 | 24 | 29 | 14 | **147** |
| 34 | 19 | 30 | 15 | 26 | 23 | **147** |
| 194| 194| 200| 200| 194| 194 |

| 15 | 30 | 19 | 34 | 23 | 26 | **147** |
|----|----|----|----|----|----|---------|
| 18 | 33 | 16 | 25 | 20 | 35 | **147** |
| 29 | 14 | 31 | 22 | 27 | 24 | **147** |
| 32 | 17 | 28 | 13 | 36 | 21 | **147** |
| 45 | 12 | 37 | 4  | 41 | 8  | **147** |
| 48 | 3  | 46 | 7  | 38 | 5  | **147** |
| 11 | 44 | 1  | 40 | 9  | 42 | **147** |
| 2  | 47 | 10 | 43 | 6  | 39 | **147** |
| 200| 200| 188| 188| 200| 200 |

Fig.43. Quasi-magic knight tours with short rows magic on 6 x 8 board.

The number of semi-magic knight tours with long rows magic is unknown but the author conjectures that it will be around 3000. J. C. Meyrignac only searched for semi-magic tours having long rows magic with tours beginning in the corner and found 5 quasi-magic tours. He didn't look for near-magic tours. The author has discovered 50 quasi-magic (48 open and 2 closed) tours with long rows magic and two examples are shown in Figure 44.

| 5 | 46 | 21 | 28 | 1 | 44 | 145 |
|---|---|---|---|---|---|---|
| 20 | 25 | 4 | 45 | 22 | 29 | 145 |
| 47 | 6 | 27 | 24 | 43 | 2 | 149 |
| 26 | 19 | 48 | 3 | 30 | 23 | 149 |
| 7 | 36 | 9 | 40 | 15 | 42 | 149 |
| 18 | 39 | 16 | 33 | 12 | 31 | 149 |
| 35 | 8 | 37 | 10 | 41 | 14 | 145 |
| 38 | 17 | 34 | 13 | 32 | 11 | 145 |
| **196** | **196** | **196** | **196** | **196** | **196** | |

| 15 | 30 | 45 | 10 | 13 | 28 | 141 |
|---|---|---|---|---|---|---|
| 46 | 11 | 14 | 29 | 44 | 9 | 153 |
| 31 | 16 | 43 | 12 | 27 | 24 | 153 |
| 2 | 47 | 18 | 25 | 8 | 41 | 141 |
| 17 | 32 | 1 | 42 | 23 | 26 | 141 |
| 48 | 3 | 36 | 19 | 40 | 7 | 153 |
| 33 | 20 | 5 | 38 | 35 | 22 | 153 |
| 4 | 37 | 34 | 21 | 6 | 39 | 141 |
| **196** | **196** | **196** | **196** | **196** | **196** | |

Fig.44. Quasi-magic knight tours with long rows magic on 6 x 8 board.

There are 80 near-magic (64 open and 16 closed) tours with long rows magic and two examples are shown in Figure 45.

| 3 | 44 | 39 | 10 | 5 | 46 | **147** |
|---|---|---|---|---|---|---|
| 38 | 9 | 4 | 45 | 40 | 11 | **147** |
| 43 | 2 | 37 | 12 | 47 | 6 | **147** |
| 36 | 19 | 8 | 41 | 30 | 13 | **147** |
| 1 | 42 | 29 | 20 | 7 | 48 | **147** |
| 18 | 35 | 22 | 25 | 14 | 31 | 145 |
| 23 | 28 | 33 | 16 | 21 | 26 | **147** |
| 34 | 17 | 24 | 27 | 32 | 15 | 149 |
| **196** | **196** | **196** | **196** | **196** | **196** | |

| 11 | 40 | 21 | 28 | 9 | 38 | **147** |
|---|---|---|---|---|---|---|
| 22 | 29 | 10 | 39 | 20 | 25 | 145 |
| 41 | 12 | 27 | 24 | 37 | 8 | 149 |
| 30 | 23 | 2 | 47 | 26 | 19 | **147** |
| 3 | 42 | 13 | 36 | 7 | 48 | 149 |
| 14 | 31 | 46 | 1 | 18 | 35 | 145 |
| 43 | 4 | 33 | 16 | 45 | 6 | **147** |
| 32 | 15 | 44 | 5 | 34 | 17 | **147** |
| **196** | **196** | **196** | **196** | **196** | **196** | |

Fig.45. Near-magic knight tours with long rows magic on 6 x 8 board.

**6 x 9 board**: The number of tours, both closed and open, are unknown. The author has enumerated 9994 arithmetically distinct semi-magic tours with short (6-cell) rows magic and two such tours are shown in Figure 46. Discerning readers must have noted that the sums of two pairs of long (9-cell) rows are identical. There are 46 such tours. There are no semi-magic tours with all the three odd (or even) long rows having the same sum. There is no quasi-magic tour.

| 1 | 36 | 7 | 48 | 19 | 54 | **165** |
|---|---|---|---|---|---|---|
| 6 | 47 | 34 | 21 | 8 | 49 | **165** |
| 35 | 2 | 37 | 18 | 53 | 20 | **165** |
| 46 | 5 | 22 | 33 | 50 | 9 | **165** |
| 23 | 38 | 3 | 52 | 17 | 32 | **165** |
| 4 | 45 | 24 | 31 | 10 | 51 | **165** |
| 39 | 30 | 41 | 26 | 13 | 16 | **165** |
| 44 | 25 | 28 | 15 | 42 | 11 | **165** |
| 29 | 40 | 43 | 12 | 27 | 14 | **165** |
| 227 | 268 | 239 | 256 | 239 | 256 | |

| 11 | 54 | 9 | 46 | 43 | 2 | **165** |
|---|---|---|---|---|---|---|
| 8 | 47 | 12 | 1 | 52 | 45 | **165** |
| 13 | 10 | 53 | 44 | 3 | 42 | **165** |
| 48 | 7 | 4 | 33 | 22 | 51 | **165** |
| 5 | 14 | 23 | 50 | 41 | 32 | **165** |
| 24 | 49 | 6 | 31 | 34 | 21 | **165** |
| 15 | 30 | 17 | 26 | 37 | 40 | **165** |
| 18 | 25 | 28 | 39 | 20 | 35 | **165** |
| 29 | 16 | 19 | 36 | 27 | 38 | **165** |
| 171 | 252 | 171 | 306 | 279 | 306 | |

Fig.46. Semi-magic knight tours on 6 x 9 board.

**6 x 10 board**: The total number of tours is unknown and tours having semi-magic properties haven't been looked into earlier. The author has enumerated 93718 arithmetically distinct semi-magic tours with short (6-cell) rows magic and two such tours are shown in Figure 47. The four non-magic long (10-cell) rows in the open tour are just +- 2 from the magic constant 305. There are 34 quasi-magic (32 open and 2 closed) tours with short rows magic and two such tours are shown in Figure 48. Discerning readers must have noted that the long rows of the open tour are +- 6 from the magic constant 305. There are 28 near-magic (12 open and 16 closed) tours with short rows magic and two such tours are shown in Figure 49. The two non-magic long rows in the open tour are just +- 6 from the magic constant 305. The number of semi-magic tours with long rows magic is unknown and the author conjectures that it will be around 150000. Figure 50 shows two near-magic tours with long rows magic. The non-magic short rows are +- 6 from the magic constant 183. These two tours are almost identical since only four digits, namely, 7, 33, and 8, 34 have swapped places. They have been highlighted in underlined italics. Such 'twin semi-magic tours' are very rare and have an aesthetic appeal.

| 17 | 12 | 9 | 52 | 49 | 44 | **183** |
|---|---|---|---|---|---|---|
| 10 | 53 | 18 | 43 | 8 | 51 | **183** |
| 13 | 16 | 11 | 50 | 45 | 48 | **183** |
| 54 | 19 | 46 | 15 | 42 | 7 | **183** |
| 23 | 14 | 21 | 40 | 47 | 38 | **183** |
| 20 | 55 | 24 | 37 | 6 | 41 | **183** |
| 25 | 22 | 27 | 34 | 39 | 36 | **183** |
| 28 | 59 | 56 | 3 | 32 | 5 | **183** |
| 57 | 26 | 33 | 30 | 35 | 2 | **183** |
| 60 | 29 | 58 | 1 | 4 | 31 | **183** |
| 307 | **305** | 303 | **305** | 307 | 303 | |

| 7 | 52 | 39 | 22 | 9 | 54 | **183** |
|---|---|---|---|---|---|---|
| 20 | 23 | 8 | 53 | 38 | 41 | **183** |
| 51 | 6 | 21 | 40 | 55 | 10 | **183** |
| 24 | 19 | 12 | 49 | 42 | 37 | **183** |
| 5 | 50 | 17 | 44 | 11 | 56 | **183** |
| 18 | 25 | 48 | 13 | 36 | 43 | **183** |
| 47 | 4 | 45 | 16 | 57 | 14 | **183** |
| 26 | 1 | 58 | 33 | 30 | 35 | **183** |
| 59 | 46 | 3 | 28 | 15 | 32 | **183** |
| 2 | 27 | 60 | 31 | 34 | 29 | **183** |
| 259 | 253 | 311 | 329 | 327 | 351 | |

Fig.47. Semi-magic knight tours with short rows magic on 6 x 10 board.

| 1 | 60 | 21 | 40 | 5 | 56 | **183** |
|---|---|---|---|---|---|---|
| 20 | 41 | 4 | 57 | 38 | 23 | **183** |
| 59 | 2 | 39 | 22 | 55 | 6 | **183** |
| 42 | 19 | 58 | 3 | 24 | 37 | **183** |
| 17 | 28 | 25 | 54 | 7 | 52 | **183** |
| 26 | 43 | 18 | 51 | 36 | 9 | **183** |
| 29 | 16 | 27 | 8 | 53 | 50 | **183** |
| 44 | 47 | 14 | 33 | 10 | 35 | **183** |
| 15 | 30 | 45 | 12 | 49 | 32 | **183** |
| 46 | 13 | 48 | 31 | 34 | 11 | **183** |
| 299 | 299 | 299 | 311 | 311 | 311 | |

| 35 | 60 | 57 | 4 | 25 | 2 | **183** |
|---|---|---|---|---|---|---|
| 58 | 5 | 36 | 1 | 56 | 27 | **183** |
| 37 | 34 | 59 | 26 | 3 | 24 | **183** |
| 6 | 9 | 32 | 53 | 28 | 55 | **183** |
| 33 | 38 | 7 | 30 | 23 | 52 | **183** |
| 8 | 31 | 10 | 51 | 54 | 29 | **183** |
| 39 | 12 | 41 | 20 | 49 | 22 | **183** |
| 42 | 15 | 50 | 11 | 46 | 19 | **183** |
| 13 | 40 | 17 | 44 | 21 | 48 | **183** |
| 16 | 43 | 14 | 47 | 18 | 45 | **183** |
| 287 | 287 | 323 | 287 | 323 | 323 | |

Fig.48. Quasi-magic knight tours with short rows magic on 6 x 10 board.

| 7 | 52 | 39 | 22 | 9 | 54 | **183** |
|---|---|---|---|---|---|---|
| 40 | 23 | 8 | 53 | 38 | 21 | **183** |
| 51 | 6 | 19 | 42 | 55 | 10 | **183** |
| 24 | 41 | 12 | 49 | 20 | 37 | **183** |
| 5 | 50 | 43 | 18 | 11 | 56 | **183** |
| 44 | 25 | 48 | 13 | 36 | 17 | **183** |
| 47 | 4 | 45 | 16 | 57 | 14 | **183** |
| 26 | 1 | 58 | 33 | 30 | 35 | **183** |
| 59 | 46 | 3 | 28 | 15 | 32 | **183** |
| 2 | 27 | 60 | 31 | 34 | 29 | **183** |
| **305** | 275 | 335 | **305** | **305** | **305** | |

| 29 | 52 | 49 | 12 | 31 | 10 | **183** |
|---|---|---|---|---|---|---|
| 50 | 13 | 30 | 9 | 48 | 33 | **183** |
| 53 | 28 | 51 | 32 | 11 | 8 | **183** |
| 14 | 17 | 26 | 45 | 34 | 47 | **183** |
| 27 | 54 | 15 | 36 | 7 | 44 | **183** |
| 16 | 25 | 18 | 43 | 46 | 35 | **183** |
| 55 | 2 | 37 | 24 | 59 | 6 | **183** |
| 38 | 19 | 60 | 1 | 42 | 23 | **183** |
| 3 | 56 | 21 | 40 | 5 | 58 | **183** |
| 20 | 39 | 4 | 57 | 22 | 41 | **183** |
| **305** | **305** | 311 | 299 | **305** | **305** | |

Fig.49. Near-magic knight tours with short rows magic on 6 x 10 board.

| 41 | 22 | 57 | 4 | 39 | 20 | **183** |
|---|---|---|---|---|---|---|
| 56 | 3 | 40 | 21 | 58 | 5 | **183** |
| 23 | 42 | 1 | 60 | 19 | 38 | **183** |
| 2 | 55 | 24 | 37 | 6 | 59 | **183** |
| 25 | 52 | 43 | 16 | 35 | 18 | 189 |
| 54 | 15 | 36 | _7_ | 44 | _33_ | 189 |
| 51 | 26 | 53 | _34_ | 17 | _8_ | 189 |
| 14 | 11 | 28 | 47 | 32 | 45 | 177 |
| 27 | 50 | 13 | 30 | 9 | 48 | 177 |
| 12 | 29 | 10 | 49 | 46 | 31 | 177 |
| **305** | **305** | **305** | **305** | **305** | **305** | |

| 41 | 22 | 57 | 4 | 39 | 20 | **183** |
|---|---|---|---|---|---|---|
| 56 | 3 | 40 | 21 | 58 | 5 | **183** |
| 23 | 42 | 1 | 60 | 19 | 38 | **183** |
| 2 | 55 | 24 | 37 | 6 | 59 | **183** |
| 25 | 52 | 43 | 16 | 35 | 18 | 189 |
| 54 | 15 | 36 | _33_ | 44 | _7_ | 189 |
| 51 | 26 | 53 | _8_ | 17 | _34_ | 189 |
| 14 | 11 | 28 | 47 | 32 | 45 | 177 |
| 27 | 50 | 13 | 30 | 9 | 48 | 177 |
| 12 | 29 | 10 | 49 | 46 | 31 | 177 |
| **305** | **305** | **305** | **305** | **305** | **305** | |

Fig.50. Near-magic knight tours with long rows magic on 6 x 10 board.

**6 x 11 board**: The total number of tours is unknown. The author has enumerated 660282 arithmetically distinct semi-magic tours with short (6-cell) rows magic and two such tours are shown in Figure 51. These two (and their reverse) are the only quasi-magic tours.

| 21 | 18 | 53 | 14 | 49 | 46 | **201** |
|---|---|---|---|---|---|---|
| 52 | 13 | 20 | 47 | 54 | 15 | **201** |
| 19 | 22 | 17 | 50 | 45 | 48 | **201** |
| 12 | 51 | 58 | 9 | 16 | 55 | **201** |
| 23 | 8 | 11 | 56 | 59 | 44 | **201** |
| 26 | 57 | 24 | 43 | 10 | 41 | **201** |
| 7 | 62 | 27 | 40 | 5 | 60 | **201** |
| 28 | 25 | 6 | 61 | 42 | 39 | **201** |
| 63 | 36 | 65 | 32 | 1 | 4 | **201** |
| 66 | 29 | 34 | 3 | 38 | 31 | **201** |
| 35 | 64 | 37 | 30 | 33 | 2 | **201** |
| 352 | 385 | 352 | 385 | 352 | 385 | |

| 33 | 4 | 35 | 62 | 65 | 2 | **201** |
|---|---|---|---|---|---|---|
| 6 | 63 | 32 | 3 | 36 | 61 | **201** |
| 31 | 34 | 5 | 64 | 1 | 66 | **201** |
| 54 | 7 | 30 | 37 | 60 | 13 | **201** |
| 9 | 38 | 53 | 14 | 29 | 58 | **201** |
| 52 | 55 | 8 | 59 | 12 | 15 | **201** |
| 39 | 10 | 51 | 16 | 57 | 28 | **201** |
| 50 | 17 | 56 | 11 | 46 | 21 | **201** |
| 43 | 40 | 45 | 22 | 27 | 24 | **201** |
| 18 | 49 | 42 | 25 | 20 | 47 | **201** |
| 41 | 44 | 19 | 48 | 23 | 26 | **201** |
| 376 | 361 | 376 | 361 | 376 | 361 | |

Fig.51. Quasi-magic knight tours with long rows magic on 6 x 11 board.

# Two-knight Magic Tours

The author has shown that a magic tour of knight is not possible up to 6 x 11 board. Before proceeding to 6 x 12 (and larger) board, it is better to consider possibilities of two-knight magic tours on smaller boards. Instead of a single knight covering the entire board, two-knight tour, as the name envisages, consists of two sequences of knight moves joined by rook moves. "The unorthodox chess-piece that combines the moves of knight and rook is known in Variant Chess as an empress and when the rook moves are restricted to single step (wazir) moves the piece (knight + wazir) is called an emperor." Jelliss [26] asserts that "there are only two two-knight emperor magic tours on the 4 x 6 board." The author has discovered that there are three, and not two, emperor magic tours on 4 x 6 board as shown in Figure 52. Jelliss [26] further asserts that "... there is one two-knight emperor tour, summing to 75 in the 6-cell lines". The author has found there are six two-knight emperor tours on 4 x 6 board as shown in Figure 53. Enumeration of tours is tricky and one can be off the mark even on smaller boards.

| 1  | 18 | 11 | 20 | **50** |
|----|----|----|----|--------|
| 10 | 21 | 2  | 17 | **50** |
| 3  | 16 | 19 | 12 | **50** |
| 22 | 9  | 6  | 13 | **50** |
| 15 | 4  | 23 | 8  | **50** |
| 24 | 7  | 14 | 5  | **50** |
| **75** | **75** | **75** | **75** | |

| 5  | 8  | 23 | 14 | **50** |
|----|----|----|----|--------|
| 24 | 15 | 4  | 7  | **50** |
| 9  | 6  | 13 | 22 | **50** |
| 16 | 19 | 12 | 3  | **50** |
| 1  | 10 | 21 | 18 | **50** |
| 20 | 17 | 2  | 11 | **50** |
| **75** | **75** | **75** | **75** | |
(Jelliss)

| 17 | 20 | 11 | 2  | **50** |
|----|----|----|----|--------|
| 10 | 1  | 18 | 21 | **50** |
| 19 | 16 | 3  | 12 | **50** |
| 6  | 9  | 22 | 13 | **50** |
| 15 | 24 | 7  | 4  | **50** |
| 8  | 5  | 14 | 23 | **50** |
| **75** | **75** | **75** | **75** | |
(Jelliss)

Fig.52. Two-knight (emperor) magic tour on 4 x 6 board.

| 1  | 22 | 3  | 18 | 11 | 20 | **75** |
|----|----|----|----|----|----|--------|
| 4  | 13 | 12 | 21 | 8  | 17 | **75** |
| 23 | 2  | 15 | 6  | 19 | 10 | **75** |
| 14 | 5  | 24 | 9  | 16 | 7  | **75** |
| 42 | 42 | 54 | 54 | 54 | 54 | |

| 1  | 18 | 11 | 22 | 3  | 20 | **75** |
|----|----|----|----|----|----|--------|
| 10 | 15 | 2  | 19 | 6  | 23 | **75** |
| 17 | 12 | 13 | 8  | 21 | 4  | **75** |
| 14 | 9  | 16 | 5  | 24 | 7  | **75** |
| 42 | 54 | 42 | 54 | 54 | 54 | |

| 5  | 18 | 3  | 22 | 7  | 20 | **75** |
|----|----|----|----|----|----|--------|
| 2  | 15 | 6  | 19 | 10 | 23 | **75** |
| 17 | 4  | 13 | 12 | 21 | 8  | **75** |
| 14 | 1  | 16 | 9  | 24 | 11 | **75** |
| 38 | 38 | 38 | 62 | 62 | 62 | |
(Jelliss)

| 15 | 6  | 19 | 10 | 23 | 2  | **75** |
|----|----|----|----|----|----|--------|
| 18 | 9  | 16 | 1  | 20 | 11 | **75** |
| 5  | 14 | 7  | 22 | 3  | 24 | **75** |
| 8  | 17 | 4  | 13 | 12 | 21 | **75** |
| 46 | 46 | 46 | 46 | 58 | 58 | |

| 15 | 10 | 23 | 6  | 19 | 2  | **75** |
|----|----|----|----|----|----|--------|
| 24 | 7  | 14 | 3  | 22 | 5  | **75** |
| 11 | 16 | 9  | 20 | 1  | 18 | **75** |
| 8  | 13 | 12 | 17 | 4  | 21 | **75** |
| 58 | 46 | 58 | 46 | 46 | 46 | |

| 17 | 4  | 13 | 12 | 21 | 8  | **75** |
|----|----|----|----|----|----|--------|
| 14 | 1  | 16 | 9  | 24 | 11 | **75** |
| 5  | 18 | 3  | 22 | 7  | 20 | **75** |
| 2  | 15 | 6  | 19 | 10 | 23 | **75** |
| 38 | 38 | 38 | 62 | 62 | 62 | |

Fig.53. Two-knight (emperor) quasi-magic tour on 4 x 6 board.

Two-knight magic tours are possible on all even size boards. So, now, let us look for single knight magic tours on 6 x 12 board.

**6 x 12 board**: The total number of tours hasn't been counted. Jelliss [19] has only given two semi-magic tours, which are also quasi-magic, with short (6-cell) rows magic by "joining together two of Awani Kumar's 6 x 6 semi-magic tours, suitably chosen." Near-magic tours haven't been investigated earlier. The author has enumerated 484 quasi-magic tours, 2056 near-magic tours and 4310970 semi-magic tours with short (6-cell) rows magic. The author is not giving any example of these tours because they have eclipsed in the light of long awaited and much sought after MAGIC TOUR. There are 8 magic tours on 6 x 12 board and are shown in Figure 54.

| 1 | 68 | 35 | 40 | 5 | 70 | **219** |
|---|---|---|---|---|---|---|
| 36 | 41 | 2 | 69 | 32 | 39 | **219** |
| 67 | 34 | 37 | 4 | 71 | 6 | **219** |
| 42 | 3 | 72 | 33 | 38 | 31 | **219** |
| 17 | 66 | 21 | 52 | 7 | 56 | **219** |
| 20 | 43 | 18 | 55 | 30 | 53 | **219** |
| 65 | 16 | 51 | 22 | 57 | 8 | **219** |
| 44 | 19 | 10 | 63 | 54 | 29 | **219** |
| 15 | 64 | 23 | 50 | 9 | 58 | **219** |
| 24 | 45 | 62 | 11 | 28 | 49 | **219** |
| 61 | 14 | 47 | 26 | 59 | 12 | **219** |
| 46 | 25 | 60 | 13 | 48 | 27 | **219** |
| **438** | **438** | **438** | **438** | **438** | **438** | |

MT.1

| 1 | 70 | 35 | 40 | 5 | 68 | **219** |
|---|---|---|---|---|---|---|
| 36 | 41 | 2 | 69 | 32 | 39 | **219** |
| 71 | 34 | 37 | 4 | 67 | 6 | **219** |
| 42 | 3 | 72 | 33 | 38 | 31 | **219** |
| 63 | 10 | 43 | 30 | 7 | 66 | **219** |
| 44 | 29 | 64 | 9 | 26 | 47 | **219** |
| 11 | 62 | 27 | 46 | 65 | 8 | **219** |
| 28 | 45 | 12 | 61 | 48 | 25 | **219** |
| 13 | 60 | 49 | 24 | 57 | 16 | **219** |
| 50 | 19 | 58 | 15 | 54 | 23 | **219** |
| 59 | 14 | 21 | 52 | 17 | 56 | **219** |
| 20 | 51 | 18 | 55 | 22 | 53 | **219** |
| **438** | **438** | **438** | **438** | **438** | **438** | |

MT.2

| 3 | 68 | 33 | 38 | 5 | 72 | **219** |
|---|---|---|---|---|---|---|
| 34 | 41 | 4 | 71 | 32 | 37 | **219** |
| 67 | 2 | 69 | 36 | 39 | 6 | **219** |
| 42 | 35 | 40 | 1 | 70 | 31 | **219** |
| 17 | 66 | 21 | 52 | 7 | 56 | **219** |
| 20 | 43 | 18 | 55 | 30 | 53 | **219** |
| 65 | 16 | 51 | 22 | 57 | 8 | **219** |
| 44 | 19 | 10 | 63 | 54 | 29 | **219** |
| 15 | 64 | 23 | 50 | 9 | 58 | **219** |
| 24 | 45 | 62 | 11 | 28 | 49 | **219** |
| 61 | 14 | 47 | 26 | 59 | 12 | **219** |
| 46 | 25 | 60 | 13 | 48 | 27 | **219** |
| **438** | **438** | **438** | **438** | **438** | **438** | |

MT.3

| 5 | 68 | 33 | 38 | 3 | 72 | **219** |
|---|---|---|---|---|---|---|
| 34 | 41 | 4 | 71 | 32 | 37 | **219** |
| 67 | 6 | 69 | 36 | 39 | 2 | **219** |
| 42 | 35 | 40 | 1 | 70 | 31 | **219** |
| 7 | 66 | 43 | 30 | 63 | 10 | **219** |
| 26 | 47 | 64 | 9 | 44 | 29 | **219** |
| 65 | 8 | 27 | 46 | 11 | 62 | **219** |
| 48 | 25 | 12 | 61 | 28 | 45 | **219** |
| 57 | 16 | 49 | 24 | 13 | 60 | **219** |
| 50 | 19 | 58 | 15 | 54 | 23 | **219** |
| 17 | 56 | 21 | 52 | 59 | 14 | **219** |
| 20 | 51 | 18 | 55 | 22 | 53 | **219** |
| **438** | **438** | **438** | **438** | **438** | **438** | |

| 17 | 58 | 19 | 54 | 15 | 56 | **219** |
|----|----|----|----|----|----|---------|
| 20 | 53 | 16 | 57 | 50 | 23 | **219** |
| 59 | 18 | 51 | 22 | 55 | 14 | **219** |
| 52 | 21 | 60 | 13 | 24 | 49 | **219** |
| 1  | 36 | 71 | 38 | 61 | 12 | **219** |
| 70 | 39 | 34 | 3  | 48 | 25 | **219** |
| 35 | 2  | 37 | 72 | 11 | 62 | **219** |
| 40 | 69 | 4  | 33 | 26 | 47 | **219** |
| 67 | 6  | 27 | 46 | 63 | 10 | **219** |
| 28 | 41 | 68 | 5  | 32 | 45 | **219** |
| 7  | 66 | 43 | 30 | 9  | 64 | **219** |
| 42 | 29 | 8  | 65 | 44 | 31 | **219** |
| **438** | **438** | **438** | **438** | **438** | **438** | |

MT.5

| 17 | 58 | 19 | 54 | 15 | 56 | **219** |
|----|----|----|----|----|----|---------|
| 20 | 53 | 16 | 57 | 50 | 23 | **219** |
| 59 | 18 | 51 | 22 | 55 | 14 | **219** |
| 52 | 21 | 60 | 13 | 24 | 49 | **219** |
| 35 | 2  | 37 | 72 | 61 | 12 | **219** |
| 70 | 39 | 34 | 3  | 48 | 25 | **219** |
| 1  | 36 | 71 | 38 | 11 | 62 | **219** |
| 40 | 69 | 4  | 33 | 26 | 47 | **219** |
| 67 | 6  | 27 | 46 | 63 | 10 | **219** |
| 28 | 41 | 68 | 5  | 32 | 45 | **219** |
| 7  | 66 | 43 | 30 | 9  | 64 | **219** |
| 42 | 29 | 8  | 65 | 44 | 31 | **219** |
| **438** | **438** | **438** | **438** | **438** | **438** | |

MT.6

| 17 | 58 | 19 | 54 | 15 | 56 | **219** |
|----|----|----|----|----|----|---------|
| 50 | 23 | 16 | 57 | 20 | 53 | **219** |
| 59 | 18 | 51 | 22 | 55 | 14 | **219** |
| 24 | 49 | 60 | 13 | 52 | 21 | **219** |
| 61 | 12 | 1  | 36 | 71 | 38 | **219** |
| 48 | 25 | 70 | 39 | 34 | 3  | **219** |
| 11 | 62 | 35 | 2  | 37 | 72 | **219** |
| 26 | 47 | 40 | 69 | 4  | 33 | **219** |
| 63 | 10 | 27 | 46 | 67 | 6  | **219** |
| 28 | 41 | 68 | 5  | 32 | 45 | **219** |
| 9  | 64 | 43 | 30 | 7  | 66 | **219** |
| 42 | 29 | 8  | 65 | 44 | 31 | **219** |
| **438** | **438** | **438** | **438** | **438** | **438** | |

MT.7

| 17 | 58 | 19 | 54 | 15 | 56 | **219** |
|----|----|----|----|----|----|---------|
| 50 | 23 | 16 | 57 | 20 | 53 | **219** |
| 59 | 18 | 51 | 22 | 55 | 14 | **219** |
| 24 | 49 | 60 | 13 | 52 | 21 | **219** |
| 61 | 12 | 35 | 2  | 37 | 72 | **219** |
| 48 | 25 | 70 | 39 | 34 | 3  | **219** |
| 11 | 62 | 1  | 36 | 71 | 38 | **219** |
| 26 | 47 | 40 | 69 | 4  | 33 | **219** |
| 63 | 10 | 27 | 46 | 67 | 6  | **219** |
| 28 | 41 | 68 | 5  | 32 | 45 | **219** |
| 9  | 64 | 43 | 30 | 7  | 66 | **219** |
| 42 | 29 | 8  | 65 | 44 | 31 | **219** |
| **438** | **438** | **438** | **438** | **438** | **438** | |

MT.8

Fig.54. Magic tours of knight on 6 x 12 board.

**6 x 16 board**: Since magic tours are not possible on an odd side (6 x 13 and 6 x 15) board and on singly-even sides (6 x 14) board, the author has ignored them and looked for magic tours on 6 x 16 board. The author has constructed over 200 magic tours on 6 x 16 board and conjectures that their total number will be in few thousands. Two magic tours are shown in Figure 55.

| 1 | 48 | 91 | 54 | 87 | 56 | 89 | 60 | 13 | 26 | 17 | 34 | 23 | 68 | 77 | 32 | **776** |
|---|----|----|----|----|----|----|----|----|----|----|----|----|----|----|----|---------|
| 92 | 53 | 46 | 3 | 90 | 59 | 86 | 57 | 16 | 35 | 14 | 25 | 78 | 33 | 22 | 67 | **776** |
| 47 | 2 | 49 | 96 | 55 | 88 | 39 | 12 | 61 | 82 | 27 | 18 | 69 | 24 | 31 | 76 | **776** |
| 52 | 93 | 4 | 45 | 42 | 9 | 58 | 85 | 36 | 15 | 70 | 79 | 28 | 73 | 66 | 21 | **776** |
| 5 | 44 | 95 | 50 | 7 | 38 | 11 | 40 | 81 | 62 | 83 | 72 | 19 | 64 | 75 | 30 | **776** |
| 94 | 51 | 6 | 43 | 10 | 41 | 8 | 37 | 84 | 71 | 80 | 63 | 74 | 29 | 20 | 65 | **776** |
| **291** | **291** | **291** | **291** | **291** | **291** | **291** | **291** | **291** | **291** | **291** | **291** | **291** | **291** | **291** | **291** | |

MT.1

| 3 | 46 | 53 | 92 | 1 | 48 | 95 | 50 | 13 | 68 | 35 | 32 | 17 | 72 | 77 | 74 | **776** |
|---|---|---|---|---|---|---|---|---|---|---|---|---|---|---|---|---|
| 54 | 91 | 2 | 47 | 94 | 51 | 12 | 37 | 34 | 31 | 14 | 69 | 78 | 75 | 16 | 71 | **776** |
| 45 | 4 | 93 | 52 | 11 | 38 | 49 | 96 | 67 | 36 | 33 | 18 | 15 | 70 | 73 | 76 | **776** |
| 90 | 55 | 44 | 5 | 58 | 87 | 10 | 39 | 30 | 61 | 64 | 79 | 82 | 27 | 24 | 21 | **776** |
| 43 | 6 | 57 | 88 | 41 | 8 | 85 | 60 | 63 | 66 | 83 | 28 | 19 | 22 | 81 | 26 | **776** |
| 56 | 89 | 42 | 7 | 86 | 59 | 40 | 9 | 84 | 29 | 62 | 65 | 80 | 25 | 20 | 23 | **776** |
| **291** | **291** | **291** | **291** | **291** | **291** | **291** | **291** | **291** | **291** | **291** | **291** | **291** | **291** | **291** | **291** | |

MT.2

Fig.55. Magic tours of knight on 6 x 16 board.

Based on above examples and since the number of tours increases rapidly with the size of the board, the author asserts that magic tours exist on all boards of size 6 x 4k for k > 2. Closed magic tours have remained elusive and readers are requested to look for them on 6 x 16 and larger boards.

**Conclusion**: Centuries old problems have been solved. Rudrata, a Kashmiri poet belonging to 9[th] century, constructed a knight tour on 4 x 8 board. Magic tour on 8 x 8 board was discovered by William Beverley in 1848 and on 12 x 12 board by Maharaja Krishnaraj Wodayer before 1868. Later, Jelliss constructed magic tour on 8 x 12 and 12 x 14 rectangular boards but a magic tour on 4 x n and 6 by n board was eluding. The author has discovered magic tours of knight on 4 x 18 and larger boards of size 4 x 20, 4 x 22 etc. According to Jelliss [23] "Magic knight's tours can be constructed on all boards 4h×4k for h and k greater than 1. (None is possible on the 4×4 or 4×8, but I'm not sure about larger cases 4×4k with k greater than 2.)" Now we are sure that magic knight's tour is possible on 4 x 4k with k = 5. According to Jelliss [26] "The question remains; are quasi-magic knight tours on boards 4m by (4n + 2) such as the 8 x 10 and 12 x 6 the best that can be achieved or are magic knight tours possible on at least some boards of these proportions?" The author has discovered 8 magic tours of knight on 6 x 12 board. The smallest rectangular board on which a magic knight tour exists is of size 4 x 18 or 6 x 12, that is, of 72 cells. It exists on all boards of size 4 x 2k for k > 8 and on 6 x 4k for k > 2. According to McGown [16] "We suspect that there is not an easy algorithm for computing the number of tours on a rectangular board, and that the problem is in fact counting hard. It would be interesting to prove or disprove this conjecture." The author has computed number of semi-magic, quasi-magic, near-magic and magic tours for various 4 x n and 6 x n board size and the results are given in Table 1 and Table 2 respectively (next page). The number of tours mentioned here are arithmetically (and not geometrically) distinct tours as the author has followed Frenicle's system of classification. Kraitchik [4] has given the number of knight's tours up to 4 x 10 board and the author has further extended it up to 4 x 14 board. Enumeration of tours, as expected, is notoriously difficult and one can be awfully off the mark even on small size boards. Readers are requested to enumerate the tours which the author couldn't. There are 1378912 semi-magic tours on 4 x 21 board. The author looked for quasi-magic tour on boards with odd side up to 4 x 21, without getting one and suspects that it doesn't exist on 4 x (2n+1) board size. Readers are also requested to prove or disprove its existence. Quasi-magic tours do exist on 6 x 11 board. Midha, [27] mentions, "Magic rectangles are well-known for their very interesting and entertaining combinatorics. Such magic rectangles have been used in designing experiments. For example, Phillips (1964, 1968a, 1968b) illustrated the use of these magic figures for the elimination of trend effects in certain classes of one-way, factorial, latin-square, and graeco-latin-square designs." Singh [28] has introduced a new approach in the field of image encryption based on knight's tour problem. Mani [29] mentions applications of magic rectangles in cryptography. Discovery of magic knight's tour on rectangular board has shown a new way to get magic rectangles and its study may lead to new applications. Ricard [30] records, "A great deal of research has gone into the knight's tour puzzle". However, a lot remains to be researched.

**Table 1**

Enumeration of Knight's Tours on 4 by n Boards

| Size | Total number of tours | Semi-magic tours | | Quasi-magic tours | | Near-magic tours | | Magic Tours |
|---|---|---|---|---|---|---|---|---|
| | | Short rows magic | Long rows magic | Short rows magic | Long rows magic | Short rows magic | Long rows magic | |
| 1 | 2 | 3 | 4 | 5 | 6 | 7 | 8 | 9 |
| 4x3 | 4 | 0 | 0 | 0 | 0 | 0 | 0 | 0 |
| 4x4 | 0 | 0 | 0 | 0 | 0 | 0 | 0 | 0 |
| 4x5 | 41 | 0 | 0 | 0 | 0 | 0 | 0 | 0 |
| 4x6 | 372 | 0 | 16 | 0 | 4 | 0 | 2 | 0 |
| 4x7 | 3189 | 9 | 0 | 0 | 0 | 0 | 0 | 0 |
| 4x8 | 15544 | 16 | 136 | 4 | 10 | 4 | 16 | 0 |
| 4x9 | 94844 | 38 | 0 | 0 | 0 | 0 | 0 | 0 |
| 4x10 | 606556 | 104 | 3102 | 12 | 30 | 16 | 22 | 0 |
| 4x11 | 3341926 | 267 | 0 | 0 | 0 | 0 | 0 | 0 |
| 4x12 | 18243164 | 608 | 58356 | 28 | 48 | 72 | 112 | 0 |
| 4x13 | 100641235 | 1444 | 0 | 0 | 0 | 0 | 0 | 0 |
| 4x14 | 526152992 | 3480 | 1092618 | 136 | 170 | 244 | 330 | 0 |
| 4x15 | ~ 2.6 billion ? | 8221 | 0 | 0 | 0 | 0 | 0 | 0 |
| 4x16 | ~ 13 billion ? | 19212 | ~ 1.6 million ? | 488 | 710 | 1012 | 1304 | 0 |
| 4x17 | ~ 66 billion ? | 45262 | 0 | 0 | 0 | 0 | 0 | 0 |
| 4x18 | ~ 340 billion ? | 213280 | ~ 26 million ? | 2624 | 2492 | 3976 | 6322 | 16 |
| 4x19 | ???? | 250247 | 0 | 0 | 0 | 0 | 0 | 0 |
| 4x20 | ???? | 587072 | ???? | 12420 | ???? | 18440 | ???? | 88 |
| 4x21 | ???? | 1378912 | ???? | 0 | 0 | 0 | 0 | 0 |
| 4x22 | ???? | 3237456 | ???? | 64120 | ???? | 81252 | ???? | 464 |
| 4x24 | ???? | 17828024 | ???? | 313584 | ???? | 388188 | ???? | 2076 |
| 4x26 | ???? | 98203312 | ???? | 1675648 | ???? | 1800728 | ???? | 9904 |
| 4x28 | ???? | 540834326 | ???? | 8536072 | ???? | 8769964 | ???? | 47456 |

**Table 2**

Enumeration of Knight's Tours on 6 by n Boards

| Size | Total number of tours | Semi-magic tours | | Quasi-magic tours | | Near-magic tours | | Magic Tours |
|---|---|---|---|---|---|---|---|---|
| | | Short rows magic | Long rows magic | Short rows magic | Long rows magic | Short rows magic | Long rows magic | |
| 1 | 2 | 3 | 4 | 5 | 6 | 7 | 8 | 9 |
| $6 \times 4$ | 372 | 0 | 16 | 0 | 4 | 0 | 2 | 0 |
| $6 \times 5$ | 9386 | 0 | 4 | 0 | 0 | 0 | 0 | 0 |
| $6 \times 6$ | 829740 | 88 | | 0 | | 0 | | 0 |
| $6 \times 7$ | ???? | 288 | 0 | 0 | 0 | 0 | 0 | 0 |
| $6 \times 8$ | ???? | 1880 | ???? | 28 | 50 | 42 | 80 | 0 |
| $6 \times 9$ | ???? | 9994 | 0 | 0 | 0 | 0 | 0 | 0 |
| $6 \times 10$ | ???? | 93718 | ???? | 34 | ???? | 28 | ???? | 0 |
| $6 \times 11$ | ???? | 660282 | 0 | 4 | 0 | 0 | 0 | 0 |
| $6 \times 12$ | ???? | 4310970 | ???? | 484 | ???? | 2056 | ???? | 8 |

**Acknowledgement**: The author is grateful to Atulit Kumar, Achint Kumar and Aushman Ojaswi for their help in preparation of this article.

**References**:

1. H. E. Dudney; *Amusements in Mathematics*, Dover 1958, pp. 101-103.

2. M. Gardner; Knights of the square table, *Mathematical Magic Show*, The Mathematical Association of America, 1989, 188–202.

3. W. W. Rouse Ball; *Mathematical Recreations and Essays*, The Macmillan Company, 1947, pp. 174-185.

4. M. Kraitchik; *Mathematical Recreations*, Dover Publications, 1953, pp. 257-266.

5. C. A. Pickover; *The Zen of Magic Squares, Circles and Stars*, Princeton University Press, 2002, pp. 210-220, 232-235.

6. M. S. Petkovic; *Famous Puzzles of Great Mathematicians*, American Mathematical Society, 2000, pp. 257-281.

7. D. Wells; *Games and Mathematics: Subtle Connections*, Cambridge University Press, 2012, pp. 76-96.

8. G.P. Jelliss; Early history of Knight's Tours, available at www.mayhematics.com.

9. H. C. von Warnsdorf (1823); *Des Rosselsprunges einfachste und allgeneinste Losung*, Schmalkalden.

10. B. R. Stonebridge (1987); The Knight's Tour of a Chessboard; *Mathematical Spectrum*, Vol.19, No.3 (1987) pp. 83-89.

11. G.P. Jelliss; Rediscovery of the Knight's Problem, available at www.mayhematics.com.

12. I. Parberry (1997); An efficient algorithm for the Knight's tour problem, *Discrete Applied Mathematics* 73 (1997) pp. 251-260.

13. O. Kyek, I. Parberry, I. Wegener (1997); Bounds on the number of knight's tours, *Discrete Applied Mathematics* 74 (1997) pp. 171-181.

14. M. M. Ismail, A. F. Z. Abidin, S. Widiyanto, M. H. Misran, M. Alice, N. A. Nordin, E. F. Shair, S. M. Mustaza, M. N. S. Zainudin (2012); Solving Knight's tour problem using Firefly algorithm, 3rd International Conference on Engineering and ICT (ICEI2012) Melaka, Malaysia, 4 - 6 April 2012.

15. P. Cull, J. D. Curtins (1978); Knight's tour revisited, *Fibonacci Quarterly* 16:3 (1978) pp. 276-286

16. K. McGown and A. Leininger (Advisor Paul Cull) (2002); Knight's Tour, Oregon State University, MIT and Oregon State University.

17. E. T. Baskoro, Subanji (2006); Methods of Constructing Knight's Tour on rectangular and diamond boards, *Math track*, Vol.2 (2006), pp. 47-55.

18. C. Flye Sainte-Marie (1877); Note sur un probleme relatif à la marche du cavalier sur l'echiquier (Séance du 18 Avril 1877), *Bulletin de la Société Mathématique de France, Année* 1876-77, vol.5, pp.144-150.

19. G.P. Jelliss; Semi-Magic Knight's Tours, available at www.mayhematics.com.

20. G.P. Jelliss ; Knight's Tours of Four-Ranks Boards, available at www.mayhematics.com.

21. G.P. Jelliss; Knight's Tours on Larger Oblong Boards, available at www.mayhematics.com.

22. G.P. Jelliss; Closed Knight's Tours of the 6 by 6 Board, available at www.mayhematics.com.

23. G.P. Jelliss; General Theory of Magic Knight's Tours, available at www.mayhematics.com.

24. A. Kumar (2002); Semi-magic knight's tour on the 6x6 board, *The Games and Puzzles Journal*, Issue # 22 (2002).

25. G.P. Jelliss; The Enumeration of Closed Knight's Tours, *The Games and Puzzles Journal, Issue* # 15 (1997).

26. G.P. Jelliss; Emperor Magic Tours, *The Games and Puzzles Journal*, Issue # 26 (2003).

27. C. K. Midha, J. P. Reyes, A. Das and L.Y. Chan (2008); On a method to construct magic rectangles of odd order, *Statistics and Applications*, Volume 6, Nos.1 & 2, 2008 (New Series), pp. 17-24.

28. M. Singh, A. Kakkar and M. Singh (2015); Image Encryption Scheme Based on Knight's Tour Problem, *Procedia Computer Science* 70 ( 2015 ) 245 – 250.

29. K. Mani and M. Viswambari (2017); A New Method of Generating Magic Rectangle, *International Journal of Computational and Applied Mathematics*, Volume 12, Number 2 (2017), pp. 391-398.

30. G. Ricard (1980); A Comparison of computer algorithms to solve for Knight's Tour, *Pi Mu Epsilon Journal*, Volume 7, No. 3, pp. 169-175.